\documentclass[11pt]{article} 
\usepackage{amssymb,amsthm,amsmath,amsfonts}
\usepackage{graphicx}
\usepackage{cite}
\usepackage{url}
  \newtheorem{theorem}{Theorem}
  
               \newtheorem{corollary}{Corollary}

               \newtheorem{remark}{Remark}
               
               \def\pf{\par\noindent {\em Proof.}~\par\noindent}
               \def\qed{~\hfill{$\square$}\pagebreak[1]\par\medskip\par}

\newcommand{\f}{{\bf{f}}}
\newcommand{\Li}{{\mbox{Lip}}}
\newcommand{\bl}{{\bf l}}

\newcommand{\La}{\mathcal{L}_{\underline{x}}}

\newcommand{\bh}{{\bf h}}
\newcommand{\bg}{{\bf g}}
\newcommand{\bj}{{\bf j}}
\newcommand{\cD}{{\cal D}}

\newcommand{\cC}{{\cal C}}

\newcommand{\R}{{\mathbb R}}

\newcommand{\I}{{\cal I}}
\newcommand{\C}{{\cal C}_{\cal L}}

\newcommand{\ux}{\underline{x}}

\newcommand{\up}{\underline{p}}
\newcommand{\ut}{\underline{t}}
\newcommand{\uy}{\underline{y}}
\newcommand{\uz}{\underline{z}}

\newcommand{\pux}{\cD_{\ux}}
\newcommand{\puy}{\cD_{\uy}}
\newcommand{\piut}{\cD_{\ut}}

\newcommand{\cS}{{\cal S}_{\cal L}}
\newcommand{\cP}{{\cal P}_{\cal L}}

\newcommand{\cT}{{\cal T}}

\begin{document}
\title{Hardy decomposition of first order Lipschitz functions by Lam\'e-Navier solutions}
\small{
\author
{Daniel Alfonso Santiesteban$^1$, Ricardo Abreu Blaya$^1$ and Daniel Alpay$^2$}
\vskip 1truecm
\date{\small$^1$ Facultad de Matem\'aticas, Universidad Aut\'onoma de Guerrero, Mexico\\$^2$ Schmid College of Science and Technology, Chapman University, USA\\Emails: danielalfonso950105@gmail.com, rabreublaya@yahoo.es, alpay@chapman.edu}

\maketitle
\begin{abstract}

\noindent
The Clifford algebra language allows us to rewrite the Lam\'e-Navier system in terms of the Euclidean Dirac operator. In this paper, the main question we shall be concerned with is whether or not a higher order Lipschitz function on the boundary $\Gamma$ of a Jordan domain $\Omega\subset\R^m$ can be decomposed into a sum of the two boundary values of a solution of the Lam\'e-Navier system with jump across $\Gamma$. Our main tool are the Hardy projections related to a singular integral operator arising in the context of Clifford analysis, which turns out to be an involution operator on the first order Lipschitz classes.
\end{abstract}

\vspace{0.3cm}

\small{

}
\section{Introduction}
The displacement vector $\vec{u}$ of the points of  a three-dimensional isotropic elastic body in the
absence of body forces is described by the Lam\'e-Navier system
\begin{equation}\label{sln}
\mu\Delta\vec{u}+(\mu+\lambda)\nabla(\nabla\cdot\vec{u})=0,
\end{equation}
where $\lambda$ and $\mu$ are the so-called Lam\'e parameters \cite{lame}. Here, the quantities $\mu > 0$ and $\lambda >-\frac{2}{3} \mu$ are elastic constants that depend on the body. The Lam\'e-Navier system (or equilibrium equation) is very important in continuum mechanics because it forms the fundamental mathematical framework for describing the behavior of elastic materials under stress. Its usefulness stems from its crucial role in linear elasticity theory and the design and analysis of structures like bridges, buildings, and mechanical components. The system allows engineers and scientists to predict structural deformation and stability, ensuring safety and efficiency in various applications (see, for example, \cite{LNS1,LNS2}).

Harmonic functions are the solutions of the second order partial differential equation $\pux\pux u=0$, where $\pux$ stands for the Dirac operator factorizing the Laplacian in $\R^m$ ($\Delta=-\pux\pux$). The solutions of the sandwich-type equation $\pux u\pux=0$, the so-called inframonogenic functions, are close connected with harmonic ones in some topics of linear elasticity theory \cite{marco1,MD,MAB1,MAB2}. These exotic functions were originally introduced in \cite{MPS2,MPS1} and more extensively studied in \cite{lavicka,wang,MAB3,MAB4,MDAS}. Indeed, the classical Lam\'e-Navier system \eqref{sln} may be reformulated as
\begin{equation}\label{lame}
\left(\frac{\mu+\lambda}{2}\right)\pux\vec{u}\pux+\left(\frac{3\mu+\lambda}{2}\right)\pux\pux\vec{u}=0,
\end{equation}
where
\[
\pux:=e_1{\frac{\partial}{\partial x_1}}+e_2{\frac{\partial}{\partial x_2}}+e_3{\frac{\partial}{\partial x_3}},
\]
 is the three-dimensional Dirac operator constructed with the standard basis of $\R^3$ that generates the Clifford algebra $\R_{0,3}$ \cite{MAB3}. A natural gene\-ralization of \eqref{lame} can be conceived if we consider the Dirac operator in $\R^m$.

The above rewritten sheds some new light on the structure of the solutions of the Lam\'e-Navier system. Indeed, every solution of \eqref{lame} can be represented by a sum of a harmonic function and an inframonogenic one, and the universal solutions of \eqref{lame} are exactly those vector fields simultaneously harmonic and inframonogenic. 

Decomposing a complex function on a closed Jordan curve $\Gamma$ as the sum of the traces of two holomorphic functions, one in each open domain defined by the curve, is a classical topic from complex analysis and represents a cornerstone in the solutions of Riemann-Hilbert problems. Indeed, such a decomposition features the Hilbert transform (a singular version of the Cauchy transform) whose $\alpha$-H\"older and $L^p$ boundedness plays a fundamental role in real harmonic analysis \cite{Da,Ga,Mu,StZy}. Such a decomposition has been extended to $m$-dimensional Euclidean spaces by means of the so-called monogenic functions (the null solutions of the Dirac operator $\pux$).

The standard higher order Lipschitz class $\Li(k+\alpha,\Gamma)$ consists of all collections of $\R$-valued continuous functions
\begin{equation}\label{LD}
\f:=\{f^{(\bj)},\,|\bj|\le k\}
\end{equation} 
defined on $\Gamma$ and satisfying the compatibility conditions
\[
\left|f^{(\bj)}(\ux)-\sum_{|\bj+\bl|\le k}\frac{f^{(\bj+\bl)}(\uy)}{\bl !}(\ux-\uy)^\bl\right|=\mathcal{O}(|\ux-\uy|^{k+\alpha-|\bj|}),\,\,\ux,\uy\in\Gamma,\,|\bj|\le k.
\]
In 1934 the American mathematician Hassler Whitney proved in \cite{Wh} that such a collection can be extended as a $C^{k,\alpha}$-smooth function on $\R^m$. For an excellent reference alone more classical lines we refer the reader to the well-known book of E. M. Stein \cite[p. 176]{St}.

As is customary, the $\R_{0,m}$-valued function classes are component-wise defined. For simplicity of notation, we continue to write $\Li(k+\alpha,\Gamma)$ for higher order Lipschitz class of $\R_{0,m}$-valued functions. For $k=0$, we have in particular $\Li(k+\alpha,\Gamma)=C^{0,\alpha}(\Gamma)$. The $\R_{0,m}$-valued class $\Li(k+\alpha,\Gamma)$ may also be defined by means of the compatibility conditions
\[
\left\|f^{(\bj)}(\ux)-\sum_{|\bl|\le k-|\bj|}\frac{f^{(\bj+\bl)}(\uy)}{\bl !}(\ux-\uy)^\bl\right\|\le  C|\ux-\uy|^{k+\alpha-|\bj|},\,\,\ux,\uy\in\Gamma,
\]
where the Clifford norm $\|\cdot\|$  is utilized and the functions $f^{\bj}$ are this time $\R_{0,m}$-valued.
If we confine ourselves to discuss the case $k=1$ our notation is considerably simplified by associating to a multi-index $\bj$ of the form $(0,0,\dots,1,\dots,0)$ the position index $j$ occupied by $1$ in $\bj$.  

A Whitney extension theorem can be stated as follows (see \cite[Theorem 3.2]{AAB} for more details)
\begin{theorem}\label{Wh}
Let $\f=\{f^{j},\,0\le j\le m\}$ be an $\R_{0,m}$-valued collection in $\Li(1+\alpha,\Gamma)$. Then, there exists a compact supported function $\tilde{f}\in C^{1,\alpha}(\R^{m},\R_{0,m})$ satisfying
\begin{itemize}
\item[(i)] $\tilde{f}|_\Gamma = f^{0},\,\partial_{x_j}\tilde{f}|_\Gamma=f^{j},\,j=1,\dots, m$,
\item[(ii)] $\tilde{f}\in C^\infty(\R^{m} \setminus \Gamma)$,
\item[(iii)] $\|\pux\tilde{f}(\ux)\| \leqslant c \, \mbox{\em dist}(\ux,\Gamma)^{\alpha-1}$, for $\ux\in\R^m\setminus\Gamma$.
\end{itemize}
\end{theorem}

In \cite{ABMM} was proved that a higher order Lipschitz function in $\Li(1+\alpha,\Gamma)$ may de decomposed into a sum of the two boundary values of a sectionally inframonogenic function in $\R^m\setminus\Gamma$. More recently, in \cite{lianet} the authors solve the problem on finding a sectionally Clifford algebra-valued harmonic function, zero at infinity and satisfying certain boundary value condition related to higher order Lipschitz functions. 

From this point of view it is quite reasonable and natural to ask whether or not a function $f\in\Li(1+\alpha,\Gamma)$ on the boundary $\Gamma$ of a Jordan domain $\Omega\subset\R^m$ can be decomposed into a sum of the two boundary values of a solution of the Lam\'e-Navier system  with jump across $\Gamma$. That is the main question we shall be concerned with in this paper. 

\section{Preliminary notions and auxiliary results}
Let us consider the $2^m$-dimensional real Clifford algebra $\R_{0,m}$ generated by $e_1,\ldots,e_m$ according to the multiplication rules $e_ie_j+e_je_i=-2\delta_{i,j},$ where $\delta_{i,j}$ is the Kronecker's symbol.  In what follows we restrict our study to the case $m>2$.

The elements of the algebra $\R_{0,m}$ have a unique representation of the form
$$a=\sum_{A\subseteq \{1,2,\ldots,m\}}a_Ae_A,$$
where $a_A\in\R$ and where we identify $e_A$ with $e_{h_1}e_{h_2}\cdots e_{h_k}$ for $A=\{h_1, h_2, \ldots,h_k\}$ $(1\leq h_1<h_2<\cdots<h_k\leq m)$ and $e_\emptyset=e_0=1$. We can embed $\R^m$ into $\R_{0,m}$ by identifying $\ux = (x_1, . . . , x_m)$ with $\ux = \sum_{i=1}^m x_ie_i$. In particular, ${Sc}[a]:=a_0$ is referred as the scalar part of $a$. Conjugation in $\R_{0,m}$ is defined as the anti-involution $a\mapsto\overline{a}$ for which $\overline{e_i}=-e_i$. A norm $\|\cdot\|$ on $\R_{0,m}$ is defined by $\|a\|^2=Sc[a\overline{a}]$ for $a\in\R_{0,m}$. 

For each $\ux\in\R^m$ it is remarkable that $\|\ux\|=|\ux|$, where $|\cdot|$ denotes the usual Euclidean norm. In general, for $a,b\in\R_{0,m}$ we have 
\begin{equation}\label{norm}
\|ab\|\le 2^{\frac{m}{2}}\|a\|\|b\|. 
\end{equation}
For an in-depth study of Clifford algebras we refer the reader to \cite{P,GM,GHS,BDS,DSS,Ryan,Ryan1,shapiro1}.

Let $\Omega\subset\R^{m}$ be a Jordan domain with sufficiently smooth boundary $\Gamma$. In what follows, we use the temporary notation $\Omega_+$ and $\Omega_-$ for the interior and exterior (containing the infinity point) domains, determined by $\Gamma$, respectively.

We will be interested in functions $f:\Omega\rightarrow$ $\R_{0,m}$, which might be written as $f(\ux)=\sum_{A}f_A(\ux)e_A$ with $f_A$ $\R$-valued. Properties, such continuity, differentiability, integrability, and so on, are ascribed coordinate-wise. In particular, we define in this way the following right modules of $\R_{0,m}$-valued functions:
\begin{itemize}
\item $C^k(\Omega,\R_{0,m}), k\in\mathbb{N}\cup\{0\}$ the right module of all $\R_{0,m}$-valued functions, $k$-times continuously differentiable in $\Omega$. We will write 
$$\partial^{\bj}:=\frac{\partial^{|\bj|}}{\partial x_1^{j_1}\partial x_2^{j_2}\dots \partial x_m^{j_m}},$$
where $\bj= (j_1,j_2,\dots , j_m)\in\left(\mathbb N\cup\{0\}\right)^{m}$ is a $m$-dimensional multi-index and 
$$|\bj|=j_1+\cdots+j_m.$$
\item $C^{k,\alpha}(\Omega,\R_{0,m})$, $\alpha \in (0, 1]$ the right module of all $\R_{0,m}$-valued functions, $k$-times $\alpha$-H\"older continuously differentiable in $\Omega$.
\end{itemize}

Let us point out that the fundamental solution (the so-called Clifford-Cauchy kernel) of the operator $\pux$ is given by
\[E_0(\ux)=\pux E_1(\ux)=-\frac{1}{\sigma_{m}}\frac{\ux}{|\ux|^{m}},\]
where
\[
E_1 (\ux)=\frac{1}{ (m-2)\sigma_m{|\ux|}^{m-2}},\,\,\ux\neq 0,
\]
is the fundamental solution of Laplacian and $\sigma_{m}$ stands for the surface area of the unit sphere in $\R^{m}$. Thus, the kernel $E_0 (\ux)$ is basic example of a two-sided monogenic function in $\R^{m}\setminus\{0\}$.

\subsection{Representation of Lam\'e-Navier solutions}
Let $f\in C^1(\Omega\cup \Gamma, \R_{0,m})$ be left monogenic in $\Omega$, then the following Clifford Cauchy type formula reproduces $f$ in $\Omega$ from its values on $\Gamma$:
\begin{equation}\label{C}
f(\ux)=\int\limits_{\Gamma}E_0(\uy-\ux)\underline{n}(\uy)f(\uy)d\uy,,\,\,\,\ux\in\Omega,
\end{equation}
where $\underline{n}(\uy)$ stands for the outward pointing unit normal at $\uy\in\Gamma$. In the case of right monogenic functions, we have the right-handed version of \eqref{C}:
\begin{equation}\label{Cr}
f(\ux)=\int\limits_{\Gamma}f(\uy)\underline{n}(\uy)E_0(\uy-\ux)d\uy,,\,\,\,\ux\in\Omega.
\end{equation}
The Clifford Cauchy kernel $E_0 (\ux)$ generates the following important integrals:
\[
({\cC}^l\,f)(\ux):=\int\limits_{\Gamma}E_0(\uy-\ux)\underline{n}(\uy)f(\uy)d\uy,\,\,\,\ux\notin\Gamma,
\]
and
\[
({\cC}^r\,f)(\ux):=\int\limits_{\Gamma}f(\uy)\underline{n}(\uy)E_0(\uy-\ux)d\uy,\,\,\,\ux\notin\Gamma,
\] 
which are reminiscent of the standard Cauchy type integral. The corresponding Teodorescu transforms are defined by:
\[
({\cT}^l\,f)(\ux):=\int\limits_{\Omega}E_0(\uy-\ux)f(\uy)d\uy,
\]
and
\[
({\cT}^r\,f)(\ux):=\int\limits_{\Omega}f(\uy)E_0(\uy-\ux)d\uy.
\] 

The Cauchy and Teodorescu transforms associated to $\pux\pux(.)$ may be defined by
\begin{equation}\label{Thar}
(\mathcal{C}_{\mathcal{H}}^lf)(\ux):=-\int_{\Gamma}E_1(\uy-\ux)\underline{n}(\uy)f(\uy)d\uy,\quad \ux\notin\Gamma,
\end{equation}
and
\begin{equation}\label{TThar}
(\mathcal{T}_{\mathcal{H}}f)(\ux):=\int_{\Omega}E_1(\uy-\ux)f(\uy)d\uy,
\end{equation}
respectively. The right-handed version of \eqref{TThar} is given by 
\begin{equation}\label{Thar2}
(\mathcal{C}_{\mathcal{H}}^rf)(\ux):=-\int_{\Gamma}f(\uy)\underline{n}(\uy)E_1(\uy-\ux)d\uy,\quad \ux\notin\Gamma.
\end{equation}
An harmonic function $f\in C^2(\Omega,\R_{0,m})\cap C^1(\Omega\cup\Gamma,\R_{0,m})$ admits in $\Omega$ the representations:
\begin{equation}
f(\ux)=(\mathcal{C}^lf)(\ux)+(\mathcal{C}_{\mathcal{H}}^l\pux f)(\ux),
\end{equation}
and
\begin{equation}
f(\ux)=(\mathcal{C}^rf)(\ux)+(\mathcal{C}_{\mathcal{H}}^r f\pux)(\ux).
\end{equation}

The Cauchy and Teodorescu transforms associated to $\pux(.)\pux$, which have been introduced in \cite{MAB1}, have the form
\begin{equation}\label{tinfra}
(\mathcal{C}_{\mathcal{I}}^rf)(\ux):=\frac{1}{2}\left(\int_{\Gamma}E_0(\uy-\ux)f(\uy)\underline{n}(\uy)(\uy-\ux)d\uy+\sum_{j=1}^me_j\int_{\Gamma} E_1(\uy-\ux)f(\uy)\underline{n}(\uy)d\uy e_j\right),
\end{equation}
and
\begin{equation}\label{ttinfra}
(\mathcal{T}_{\mathcal{I}}f)(\ux):=-\frac{1}{2}\left(\int_{\Omega}E_0(\uy-\ux)f(\uy)(\uy-\ux)d\uy+\sum_{j=1}^me_j\int_\Omega E_1(\uy-\ux)f(\uy)d\uy e_j\right).
\end{equation}
We also have a left-handed version of \eqref{tinfra}, which is
\begin{equation}\label{tinfra2}
(\mathcal{C}_{\mathcal{I}}^lf)(\ux):=\frac{1}{2}\left(\int_{\Gamma}E_0(\uy-\ux)\underline{n}(\uy)f(\uy)(\uy-\ux)d\uy+\sum_{j=1}^me_j\int_{\Gamma}E_1(\uy-\ux) \underline{n}(\uy)f(\uy)d\uy e_j\right).
\end{equation}

An inframonogenic function $f\in C^2(\Omega,\R_{0,m})\cap C^1(\Omega\cup \Gamma, \R_{0,m})$ admits in $\Omega$ the representations:  
\begin{eqnarray}\label{cfi}
f(\ux)=(\mathcal{C}^lf)(\ux)+(\mathcal{C}_{\mathcal{I}}^r\pux f)(\ux)
\end{eqnarray}
and
\begin{eqnarray}\label{cfi}
f(\ux)=(\mathcal{C}^rf)(\ux)+(\mathcal{C}_{\mathcal{I}}^l f\pux)(\ux).
\end{eqnarray}

We will denote by $\La$ the Cliffordian Lam\'e-Navier operator
$$\left(\frac{\mu+\lambda}{2}\right)\pux(.)\pux+\left(\frac{3\mu+\lambda}{2}\right)\pux\pux(.).$$
We will commonly refer to the null solutions of the operator $\La$ as Lam\'e-Navier solutions. 
It is easy to see that this operator can be factorized as
\begin{equation}
\La=\mathcal{M}_{\ux}\pux(.)=\pux\mathcal{M}_{\ux}(.)=\overline{\mathcal{M}}_{\ux}(.)\pux=(.)\pux\overline{\mathcal{M}}_{\ux},
\end{equation}
where $\mathcal{M}_{\ux}:=\left(\frac{\mu+\lambda}{2}\right)(.)\pux+\left(\frac{3\mu+\lambda}{2}\right)\pux(.)$ and $\overline{\mathcal{M}}_{\ux}:=\left(\frac{\mu+\lambda}{2}\right)\pux(.)+\left(\frac{3\mu+\lambda}{2}\right)(.)\pux$.
The Cauchy and Teodorescu transforms associated to $\La$ will be defined as a linear combination of the corresponding integral operators  \eqref{Thar}, \eqref{TThar}, \eqref{tinfra} and \eqref{ttinfra}, a linear combination which involves the Lam\'e parameters. Indeed, we define 
\begin{equation}\label{CLiz}
(\mathcal{C}_{\mathcal{L}}^l f)(\ux):=-\frac{\mu+\lambda}{2\mu(2\mu+\lambda)}(\mathcal{C}_{\mathcal{I}}^lf)(\ux)+\frac{3\mu+\lambda}{2\mu (2\mu+\lambda)}(\mathcal{C}_{\mathcal{H}}^lf)(\ux),
\end{equation}
and
\begin{equation}
(\mathcal{T}_{\mathcal{L}} f)(\ux):=-\frac{\mu+\lambda}{2\mu (2\mu+\lambda)}(\mathcal{T}_{\mathcal{I}}f)(\ux)+\frac{3\mu+\lambda}{2\mu (2\mu+\lambda)}(\mathcal{T}_{\mathcal{H}}f)(\ux).
\end{equation}
The right-handed version of \eqref{CLiz} is defined by
\begin{equation}\label{CLde}
(\mathcal{C}_{\mathcal{L}}^r f)(\ux):=-\frac{\mu+\lambda}{2\mu(2\mu+\lambda)}(\mathcal{C}_{\mathcal{I}}^rf)(\ux)+\frac{3\mu+\lambda}{2\mu (2\mu+\lambda)}(\mathcal{C}_{\mathcal{H}}^rf)(\ux).
\end{equation}

\begin{theorem}[Borel-Pompeiu formula] Let $f\in C^2(\overline{\Omega})$. Then
\begin{align}
&(\mathcal{C}^lf)(\ux)+\frac{3\mu+\lambda}{2\mu(2\mu+\lambda)}(\mathcal{C}_{\mathcal{H}}^l\mathcal{M}_{\ux}f)(\ux)-\frac{\mu+\lambda}{2\mu(2\mu+\lambda)}(\mathcal{C}_{\mathcal{I}}^r\overline{\mathcal{M}}_{\ux}f)(\ux)+\mathcal{T}_{\mathcal{L}}(\La f)(\ux)\nonumber\\
&=\left\{\begin{array}{rl}
f(\ux),&\text{if}\;\ux\in\Omega,\\
0,&\text{if}\;\ux\in\R^m\setminus\overline{\Omega}.
\end{array}\right.\label{bpform}
\end{align}
\end{theorem}
\pf
We prove only the statement for $\ux\in\Omega$. The proof for $\ux\in\R^m\setminus\overline{\Omega}$ is completely analogous.

Using the Borel-Pompeiu formulas involving the operators $\pux\pux(.)$ and $\pux(.)\pux$ (see \cite{bory} and \cite{MAB1}):
\begin{equation}\label{fbph}
f(\ux)=(\mathcal{C}^l f)(\ux)+(\mathcal{C}_{\mathcal{H}}^l\pux f)(\ux)+(\mathcal{T}_{\mathcal{H}}\pux\pux f)(\ux) 
\end{equation}
and
\begin{equation}\label{fbpi}
f(\ux)=(\mathcal{C}^lf)(\ux)+(\mathcal{C}_\mathcal{I}^r\pux f)(\ux)+(\mathcal{T}_\mathcal{I}\pux f\pux)(\ux),
\end{equation}
we have
\begin{align*}
\mathcal{T}_{\mathcal{L}}(\La f)(\ux)&=-\frac{(\mu+\lambda)^2}{4\mu(2\mu+\lambda)}(\mathcal{T}_\mathcal{I}\pux f\pux)(\ux)-\frac{(\mu+\lambda)(3\mu+\lambda)}{4\mu(2\mu+\lambda)}(\mathcal{T}_\mathcal{I}\pux\pux f)(\ux)\\
&\quad+\frac{(3\mu+\lambda)(\mu+\lambda)}{4\mu(2\mu+\lambda)}(\mathcal{T}_{\mathcal{H}}\pux f\pux)(\ux)+\frac{(3\mu+\lambda)^2}{4\mu(2\mu+\lambda)}(\mathcal{T}_{\mathcal{H}}\pux\pux f)(\ux)\\
&=f(\ux)+\frac{(\mu+\lambda)^2}{4\mu(2\mu+\lambda)}(\mathcal{C}^lf)(\ux)-\frac{(3\mu+\lambda)^2}{4\mu(2\mu+\lambda)}(\mathcal{C}^lf)(\ux)\\
&\quad+ \frac{(\mu+\lambda)^2}{4\mu(2\mu+\lambda)}(\mathcal{C}_{\mathcal{I}}^r \pux f)(\ux)-\frac{(3\mu+\lambda)^2}{4\mu(2\mu+\lambda)}(\mathcal{C}_{\mathcal{H}}^l\pux f)(\ux)\\
&\quad+\frac{(3\mu+\lambda)(\mu+\lambda)}{4\mu(2\mu+\lambda)}[(\mathcal{T}_{\mathcal{H}}\pux f\pux)(\ux)-(\mathcal{T}_{\mathcal{I}}\pux\pux f)(\ux)]. 
\end{align*}
Now, applying the equations \eqref{fbph} and \eqref{fbpi}, and using the identity
\begin{equation}
(\mathcal{T}^l f)\pux=\pux(\mathcal{T}^r f),
\end{equation}
we get
\begin{equation}
(\mathcal{T}_{\mathcal{H}}\pux f\pux)(\ux)-(\mathcal{T}_{\mathcal{I}}\pux\pux f)(\ux)=(\mathcal{C}_{\mathcal{I}}^r f\pux)(\ux)-(\mathcal{C}_{\mathcal{H}}^l f\pux)(\ux). 
\end{equation}

Therefore, for $\ux\in\Omega$,
\begin{align*}
&f(\ux)\\&=(\mathcal{C}^lf)(\ux)-\frac{(\mu+\lambda)^2}{4\mu(2\mu+\lambda)}(\mathcal{C}_{\mathcal{I}}^r\pux f)(\ux)+\frac{(3\mu+\lambda)^2}{4\mu(2\mu+\lambda)}(\mathcal{C}_{\mathcal{H}}^l\pux f)(\ux)\\
&\quad+\frac{(3\mu+\lambda)(\mu+\lambda)}{4\mu(2\mu+\lambda)}(\mathcal{C}_{\mathcal{H}}^lf\pux)(\ux)-\frac{(3\mu+\lambda)(\mu+\lambda)}{4\mu(2\mu+\lambda)}(\mathcal{C}_{\mathcal{I}}^r f\pux)(\ux)+\mathcal{T}_{\mathcal{L}}(\La f)(\ux)\\
&=(\mathcal{C}^lf)(\ux)+\frac{3\mu+\lambda}{2\mu(2\mu+\lambda)}(\mathcal{C}_{\mathcal{H}}^l\mathcal{M}_{\ux}f)(\ux)-\frac{\mu+\lambda}{2\mu(2\mu+\lambda)}(\mathcal{C}_{\mathcal{I}}^r\overline{\mathcal{M}}_{\ux}f)(\ux)+\mathcal{T}_{\mathcal{L}}(\La f)(\ux),
\end{align*}	
and we are done.
\qed
Following the same lines as the previous proof, we can also prove the following representations
\begin{align}
&(\mathcal{C}^rf)(\ux)+\frac{3\mu+\lambda}{2\mu(2\mu+\lambda)}(\mathcal{C}_{\mathcal{H}}^r\overline{\mathcal{M}}_{\ux}f)(\ux)-\frac{\mu+\lambda}{2\mu(2\mu+\lambda)}(\mathcal{C}_{\mathcal{I}}^l\mathcal{M}_{\ux}f)(\ux)+\mathcal{T}_{\mathcal{L}}(\La f)(\ux)\nonumber\\
&=\left\{\begin{array}{rl}
f(\ux),&\text{if}\;\ux\in\Omega,\\
0,&\text{if}\;\ux\in\R^m\setminus\overline{\Omega},
\end{array}\right.\\
&\frac{(3\mu+\lambda)^2}{4\mu (2\mu+\lambda)}(\mathcal{C}^lf)(\ux)-\frac{(\mu+\lambda)^2}{4\mu (2\mu+\lambda)}(\mathcal{C}^rf)(\ux)+\mathcal{C}_{\mathcal{L}}^l\left(\mathcal{M}_{\ux}f\right)(\ux)+\mathcal{T}_{\mathcal{L}}(\La f)(\ux)\nonumber\\
&+\frac{(3\mu+\lambda)(\mu+\lambda)}{4\mu (2\mu+\lambda)}\left(\int_{\Gamma}E_0(\uy-\ux)f(\uy)\underline{n}(\uy)d\uy-\int_{\Gamma}\underline{n}(\uy)f(\uy)E_0(\uy-\ux)d\uy\right)\nonumber\\
&=\left\{\begin{array}{rl}
f(\ux),&\text{if}\;\ux\in\Omega,\\
0,&\text{if}\;\ux\in\R^m\setminus\overline{\Omega},
\end{array}\right.\label{bpform}\\
&\frac{(3\mu+\lambda)^2}{4\mu (2\mu+\lambda)}(\mathcal{C}^rf)(\ux)-\frac{(\mu+\lambda)^2}{4\mu (2\mu+\lambda)}(\mathcal{C}^lf)(\ux)+\mathcal{C}_{\mathcal{L}}^r\left(\overline{\mathcal{M}}_{\ux}f\right)(\ux)+\mathcal{T}_{\mathcal{L}}(\La f)(\ux)\nonumber\\
&+\frac{(3\mu+\lambda)(\mu+\lambda)}{4\mu (2\mu+\lambda)}\left(\int_{\Gamma}\underline{n}(\uy)f(\uy)E_0(\uy-\ux)d\uy-\int_{\Gamma}E_0(\uy-\ux)f(\uy)\underline{n}(\uy)d\uy\right)\nonumber\\
&=\left\{\begin{array}{rl}
f(\ux),&\text{if}\;\ux\in\Omega,\\
0,&\text{if}\;\ux\in\R^m\setminus\overline{\Omega}.
\end{array}\right.\label{bpform}
\end{align}
It is also easy to prove that
\begin{align}\label{relutil}
\mathcal{M}_{\ux}\left[(\mathcal{C}^lf)(\ux)+\frac{3\mu+\lambda}{2\mu(2\mu+\lambda)}(\mathcal{C}_{\mathcal{H}}^l\mathcal{M}_{\ux}f)(\ux)-\frac{\mu+\lambda}{2\mu(2\mu+\lambda)}(\mathcal{C}_{\mathcal{I}}^r\overline{\mathcal{M}}_{\ux}f)(\ux)\right]=(\mathcal{C}^l\mathcal{M}_{\ux}f)(\ux).
\end{align}
Indeed, we have

\begin{align*}
&\mathcal{M}_{\ux}\left[(\mathcal{C}^lf)(\ux)+\frac{3\mu+\lambda}{2\mu(2\mu+\lambda)}(\mathcal{C}_{\mathcal{H}}^l\mathcal{M}_{\ux}f)(\ux)-\frac{\mu+\lambda}{2\mu(2\mu+\lambda)}(\mathcal{C}_{\mathcal{I}}^r\overline{\mathcal{M}}_{\ux}f)(\ux)\right]\\
&=\mathcal{M}_{\underline{x}}\left[\frac{(3\mu+\lambda)^2}{4\mu (2\mu+\lambda)}(\mathcal{C}^lf)(\ux)-\frac{(\mu+\lambda)^2}{4\mu (2\mu+\lambda)}(\mathcal{C}^rf)(\ux)+\mathcal{C}_{\mathcal{L}}^l\left(\mathcal{M}_{\ux}f\right)(\ux)\nonumber\nonumber\right.\\
&\quad\left.
+\frac{(3\mu+\lambda)(\mu+\lambda)}{4\mu (2\mu+\lambda)}\left(\int_{\Gamma}E_0(\uy-\ux)f(\uy)\underline{n}(\uy)d\uy-\int_{\Gamma}\underline{n}(\uy)f(\uy)E_0(\uy-\ux)d\uy\right)\right]\\
&=\frac{(3\mu+\lambda)^2(\mu+\lambda)}{8\mu(2\mu+\lambda)}[(\mathcal{C}^lf)(\ux)]\pux-\frac{(\mu+\lambda)^2(3\mu+\lambda)}{8\mu(2\mu+\lambda)}\pux[(\mathcal{C}^rf)(\ux)]\\
&\quad+\frac{(3\mu+\lambda)(\mu+\lambda)^2}{8\mu(2\mu+\lambda)}\left(\int_\Gamma E_0(\uy-\ux)f(\uy)\underline{n}(\uy)d\uy\right)\pux\\&\quad-\frac{(3\mu+\lambda)^2(\mu+\lambda)}{8\mu(2\mu+\lambda)}\pux\left(\int_\Gamma \underline{n}(\uy)f(\uy)E_0(\uy-\ux)d\uy\right)\\
&\quad-\frac{(\mu+\lambda)^2}{4\mu(2\mu+\lambda)}(\mathcal{C}_\mathcal{I}^l\mathcal{M}_{\ux}f)\pux-\frac{(\mu+\lambda)(3\mu+\lambda)}{4\mu(2\mu+\lambda)}\pux(\mathcal{C}_\mathcal{I}^l\mathcal{M}_{\ux}f)\\
&\quad+\frac{(\mu+\lambda)(3\mu+\lambda)}{4\mu(2\mu+\lambda)}(\mathcal{C}_\mathcal{H}^l\mathcal{M}_{\ux}f)\pux+\frac{(3\mu+\lambda)^2}{4\mu(2\mu+\lambda)}\pux(\mathcal{C}_\mathcal{H}^l\mathcal{M}_{\ux}f)\\
&=(\mathcal{C}^l\mathcal{M}_{\ux}f)(\ux).
\end{align*}
Similarly,
\begin{align}\label{relutil2}
\overline{\mathcal{M}}_{\ux}\left[(\mathcal{C}^lf)(\ux)+\frac{3\mu+\lambda}{2\mu(2\mu+\lambda)}(\mathcal{C}_{\mathcal{H}}^l\mathcal{M}_{\ux}f)(\ux)-\frac{\mu+\lambda}{2\mu(2\mu+\lambda)}(\mathcal{C}_{\mathcal{I}}^r\overline{\mathcal{M}}_{\ux}f)(\ux)\right]=(\mathcal{C}^r\overline{\mathcal{M}}_{\ux}f)(\ux).
\end{align}

\begin{corollary}[Cauchy integral formula]
Let $f\in C^2(\Omega)\cap C^1(\overline{\Omega})$. If $\La f=0$ in $\Omega$ then
\begin{equation}\label{CLame}
f(\ux)=(\mathcal{C}^lf)(\ux)+\frac{3\mu+\lambda}{2\mu(2\mu+\lambda)}(\mathcal{C}_{\mathcal{H}}^l\mathcal{M}_{\ux}f)(\ux)-\frac{\mu+\lambda}{2\mu(2\mu+\lambda)}(\mathcal{C}_{\mathcal{I}}^r\overline{\mathcal{M}}_{\ux}f)(\ux).
\end{equation} 
Some rewritings of \eqref{CLame} are as follows
\begin{align}
f(\ux)&=(\mathcal{C}^rf)(\ux)+\frac{3\mu+\lambda}{2\mu(2\mu+\lambda)}(\mathcal{C}_{\mathcal{H}}^r\overline{\mathcal{M}}_{\ux}f)(\ux)-\frac{\mu+\lambda}{2\mu(2\mu+\lambda)}(\mathcal{C}_{\mathcal{I}}^l\mathcal{M}_{\ux}f)(\ux),\\
f(\ux)&=\frac{(3\mu+\lambda)^2}{4\mu (2\mu+\lambda)}(\mathcal{C}^lf)(\ux)-\frac{(\mu+\lambda)^2}{4\mu (2\mu+\lambda)}(\mathcal{C}^rf)(\ux)+\mathcal{C}_{\mathcal{L}}^l\left(\mathcal{M}_{\ux}f\right)(\ux)\nonumber\\
&\quad+\frac{(3\mu+\lambda)(\mu+\lambda)}{4\mu (2\mu+\lambda)}\left(\int_{\Gamma}E_0(\uy-\ux)f(\uy)\underline{n}(\uy)d\uy-\int_{\Gamma}\underline{n}(\uy)f(\uy)E_0(\uy-\ux)d\uy\right),\\
f(\ux)&=\frac{(3\mu+\lambda)^2}{4\mu (2\mu+\lambda)}(\mathcal{C}^rf)(\ux)-\frac{(\mu+\lambda)^2}{4\mu (2\mu+\lambda)}(\mathcal{C}^lf)(\ux)+\mathcal{C}_{\mathcal{L}}^r\left(\overline{\mathcal{M}}_{\ux}f\right)(\ux)\nonumber\\
&\quad+\frac{(3\mu+\lambda)(\mu+\lambda)}{4\mu (2\mu+\lambda)}\left(\int_{\Gamma}\underline{n}(\uy)f(\uy)E_0(\uy-\ux)d\uy-\int_{\Gamma}E_0(\uy-\ux)f(\uy)\underline{n}(\uy)d\uy\right).
\end{align}
\end{corollary}

For the purpose of this paper the following representation formula of Lam\'e-Navier solutions in the exterior domain $\Omega_-$ will be needed.
\begin{theorem}\label{formulaEx}
Let $f\in C^2(\Omega_-, \R_{0,m})\cap C^1(\Omega_{-}\cup \Gamma, \R_{0,m})$ be a Lam\'e-Navier solution in $\Omega_-$ such that $f(\infty)$ exists and 
\[
\|\pux f(\ux)\|=o\bigg(\frac{1}{|\ux|}\bigg)\,\,{\mbox{as}}\,\,\ux\to\infty.
\]
Then for $\ux\in\Omega_-$ we have
\begin{eqnarray}\label{cfe}
f(\ux)=-(\mathcal{C}^lf)(\ux)-\frac{3\mu+\lambda}{2\mu(2\mu+\lambda)}(\mathcal{C}_{\mathcal{H}}\mathcal{M}_{\ux}f)(\ux)+\frac{\mu+\lambda}{2\mu(2\mu+\lambda)}(\mathcal{C}_{\mathcal{I}}\overline{\mathcal{M}}_{\ux}f)(\ux)+f(\infty).
\end{eqnarray}
\end{theorem}
\pf Let $\ux\in\Omega_-$ and choose a sphere $C_R(\ux):=\{\uy\in\R^m:\,|\uy-\ux|\le R\}$ with sufficiently large radius $R$ which
contains the set $\Omega\cup \Gamma$.

The Cauchy integral formula \eqref{CLame} applied to $\Omega_-\cap\{\uy\in\R^m:\,|\uy-\ux|< R\}$ yields
\begin{align*}
&f(\ux)=-\int_\Gamma E_0(\uy-\ux)\underline{n}(\uy)f(\uy)d\uy+\frac{3\mu+\lambda}{2\mu(2\mu+\lambda)}\int_\Gamma E_1(\uy-\ux)\underline{n}(\uy)[\mathcal{M}_{\uy}f(\uy)]d\uy\\
&+\frac{\mu+\lambda}{4\mu(2\mu+\lambda)}\left(\int_{\Gamma}E_0(\uy-\ux)[\overline{\mathcal{M}}_{\uy}f(\uy)]\underline{n}(\uy)(\uy-\ux)d\uy+\sum_{j=1}^me_j\int_{\Gamma} E_1(\uy-\ux)[\overline{\mathcal{M}}_{\uy}f(\uy)]\underline{n}(\uy)d\uy e_j\right)\\
&+\int_{\partial C_R} E_0(\uy-\ux)\underline{n}(\uy)f(\uy)d\uy-\frac{3\mu+\lambda}{2\mu(2\mu+\lambda)}\int_{\partial C_R} E_1(\uy-\ux)\underline{n}(\uy)[\mathcal{M}_{\uy}f(\uy)]d\uy\\
&-\frac{\mu+\lambda}{4\mu(2\mu+\lambda)}\left(\int_{\partial C_R}E_0(\uy-\ux)[\overline{\mathcal{M}}_{\uy}f(\uy)]\underline{n}(\uy)(\uy-\ux)d\uy+\sum_{j=1}^me_j\int_{\partial C_R} E_1(\uy-\ux)[\overline{\mathcal{M}}_{\uy}f(\uy)]\underline{n}(\uy)d\uy e_j\right).
\end{align*}
From our assumptions we have that
\[
\int\limits_{\partial C_R}E_0(\uy-\ux)\underline{n}(\uy)f(\uy)d\uy \rightarrow f(\infty),
\]
while the remaining integrals on $\partial C_R$ vanish as $R\to\infty$ and the assertion is proved.\qed  
\subsection{The Cauchy type integral for Lam\'e-Navier solutions}
The previously introduced Lipschitz class $\Li(1+\alpha,\Gamma)$ is quite appropriate to define on it a sort of Cauchy type integral arising from \eqref{cfi}. More precisely, given $\f\in\Li(1+\alpha,\Gamma)$, the Cauchy type integral of $\f$ might be defined as
\begin{align}\label{CT}
&\C^0[\f](\ux)=\int\limits_\Gamma E_0(\uy-\ux)\underline{n}(\uy)f^0(\uy)d\uy-\frac{3\mu+\lambda}{2\mu(2\mu+\lambda)}\int\limits_\Gamma E_1(\uy-\ux)\underline{n}(\uy)[\mathcal{M}_{\uy}\tilde{f}(\uy)]d\uy\nonumber\\
&-\frac{\mu+\lambda}{4\mu(2\mu+\lambda)}\left(\int\limits_\Gamma E_0(\uy-\ux)[\overline{\mathcal{M}}_{\uy}\tilde{f}(\uy)]\underline{n}(\uy)(\uy-\ux)d\uy+\sum_{j=1}^me_j\int\limits_\Gamma E_1(\uy-\ux)[\overline{\mathcal{M}}_{\uy}\tilde{f}(\uy)]\underline{n}(\uy)d\uy e_j\right).
\end{align}
Here and subsequently, $\tilde{f}$ denotes the Whitney extension of $\f$, following Theorem \ref{Wh}. The restrictions of $\mathcal{M}_{\uy}\tilde{f}$ and $\overline{\mathcal{M}}_{\uy}\tilde{f}$ to $\Gamma$ give
\[
\mathcal{M}_{\uy}f(\uy)=\left(\frac{\mu+\lambda}{2}\right)\sum_{j=1}^m f^{j}(\uy)e_j+\left(\frac{3\mu+\lambda}{2}\right)\sum_{j=1}^m e_jf^{j}(\uy),
\]
and
\[
\overline{\mathcal{M}}_{\uy}\tilde{f}(\uy)=\left(\frac{3\mu+\lambda}{2}\right)\sum_{j=1}^m f^{j}(\uy)e_j+\left(\frac{\mu+\lambda}{2}\right)\sum_{j=1}^m e_jf^{j}(\uy).
\]

Applying \eqref{relutil} we obtain 
\begin{equation}\label{i-l}
\mathcal{M}_{\ux}[\C^0[\f](\ux)]={\cC}^l[\mathcal{M}_{\ux}\tilde{f}](\ux),\,\,\ux\in\R^m\setminus\Gamma,
\end{equation}
It is a simple matter to see that $\C^0[\f](\ux)$ is a Lam\'e-Navier solution in $\R^m\setminus\Gamma$, i.e. $\mathcal{L}_{\ux}[\C^0[\f](\ux)]=0$, $\ux\in\R^m\setminus\Gamma$. Moreover, the remaining integrals in \eqref{CT} (with the exception of the first) become weakly singular as $\ux$ approaches to $\Gamma$, which is easy to check. Consequently, all these integrals do not have jump discontinuity in $\Gamma$. Combining this fact with the classical Plemelj-Sokhotski formulas applied to the first Cliffordian-Cauchy type integral in \eqref{CT} we conclude that
\begin{equation}\label{Jump}
[\C^0\f]^+(\ux)-[\C^0\f]^-(\ux)=\tilde{f}(\ux)=f^{0}(\ux),\,\,\,\ux\in\Gamma,
\end{equation} 
where 
\[
[\C^0\f]^\pm(\ux):=\lim_{\Omega_\pm\ni\uz\to\ux}[\C^0\f](\uz).
\]

The singular version of $[\C^0\f]$ on $\Gamma$, denoted by $\cS^0\f$, is given by
\begin{equation}\label{Jump2}
\cS^0\f(\ux)=[\C^0\f]^+(\ux)+[\C^0\f]^-(\ux),\,\,\ux\in\Gamma,
\end{equation}
or equivalently
\begin{align}\label{S0}
&\cS^0\f(\ux)=2\int\limits_\Gamma E_0(\uy-\ux)\underline{n}(\uy)f^0(\uy)d\uy-\frac{3\mu+\lambda}{\mu(2\mu+\lambda)}\int\limits_\Gamma E_1(\uy-\ux)\underline{n}(\uy)[\mathcal{M}_{\uy}\tilde{f}(\uy)]d\uy\nonumber\\
&-\frac{\mu+\lambda}{2\mu(2\mu+\lambda)}\left(\int\limits_\Gamma E_0(\uy-\ux)[\overline{\mathcal{M}}_{\uy}\tilde{f}(\uy)]\underline{n}(\uy)(\uy-\ux)d\uy+\sum_{j=1}^me_j\int\limits_\Gamma E_1(\uy-\ux)[\overline{\mathcal{M}}_{\uy}\tilde{f}(\uy)]\underline{n}(\uy)d\uy e_j\right).
\end{align}
where the integrals in \eqref{S0} are taken in the sense of the Cauchy principal value. For brevity, the notation $\text{p.v.}$ before these integrals is omitted.

Let $\f\in\Li(1+\alpha,\Gamma)$, then for $\uy\in\R^m,\,\ux\in\Gamma$ we have (see \cite{St} for more details)
\begin{equation}\label{poly}
\tilde{f}(\uy)=P_{\ux}(\uy)+R_{\ux}(\uy),
\end{equation}
where
\[
P_{\ux}(\uy)=f^0(\ux)+\sum\limits_{j=1}^m f^j(\ux)(y_j-x_j).
\]
Since $P_{\ux}(\uy)$ is a polynomial in $\R^m$ such that $\mathcal{L}_{\uy} P_{\ux}(\uy)=0$, the combination of formula \eqref{CLame} with the obvious equality $P_{\ux}(\ux)=f^0(\ux)$ allows us to rephrase \eqref{S0} as follows
\begin{align*}
&\cS^0\f(\ux)=f^0(\ux)+2\int\limits_\Gamma E_0(\uy-\ux)\underline{n}(\uy)R_{\ux}(\uy)d\uy-\frac{3\mu+\lambda}{\mu(2\mu+\lambda)}\int\limits_\Gamma E_1(\uy-\ux)\underline{n}(\uy)[\mathcal{M}_{\uy}R_{\ux}(\uy)]d\uy\nonumber\\
&-\frac{\mu+\lambda}{2\mu(2\mu+\lambda)}\left(\int\limits_\Gamma E_0(\uy-\ux)[\overline{\mathcal{M}}_{\uy}R_{\ux}(\uy)]\underline{n}(\uy)(\uy-\ux)d\uy+\sum_{j=1}^me_j\int\limits_\Gamma E_1(\uy-\ux)[\overline{\mathcal{M}}_{\uy}R_{\ux}(\uy)]\underline{n}(\uy)d\uy e_j\right).
\end{align*}
Let us introduce one more piece of notation: 
\[
E_0^j(\uy-\ux)=\partial_{z_j} E_0(\uy-\uz)|_{\uz=\ux},\,\,E_1^j(\uy-\ux)=\partial_{z_j} E_1(\uy-\uz)|_{\uz=\ux} 
\]
and set
\begin{equation}\label{SI}
\cS\f=\{\cS^{j}\f:\,0\le j\le m\},
\end{equation}
so that for $j>0$ we have
\begin{align*}\label{Sj}
&\cS^{j}\f(\ux)=f^j(\ux)+2\int\limits_\Gamma E_0^j(\uy-\ux)\underline{n}(\uy)R_{\ux}(\uy)d\uy-\frac{3\mu+\lambda}{\mu(2\mu+\lambda)}\int\limits_\Gamma E_1^j(\uy-\ux)\underline{n}(\uy)[\mathcal{M}_{\uy}R_{\ux}(\uy)]d\uy\nonumber\\
&-\frac{\mu+\lambda}{2\mu(2\mu+\lambda)}\left(\int\limits_\Gamma E_0^j(\uy-\ux)[\overline{\mathcal{M}}_{\uy}R_{\ux}(\uy)]\underline{n}(\uy)(\uy-\ux)d\uy+\sum_{j=1}^me_j\int\limits_\Gamma E_1^j(\uy-\ux)[\overline{\mathcal{M}}_{\uy}R_{\ux}(\uy)]\underline{n}(\uy)d\uy e_j\right)\\
&+\frac{\mu+\lambda}{2\mu(2\mu+\lambda)}\int_\Gamma E_0(\uy-\ux)[\overline{\mathcal{M}}_{\uy}R_{\ux}(\uy)]\underline{n}(\uy)e_jd\uy.
\end{align*}
The following identities will be used frequently from now on
\begin{equation}\label{eqprim}
\frac{\mu+\lambda}{2}\sum\limits_{j=1}^m\cS^j\f(\ux)e_j+\frac{3\mu+\lambda}{2}\sum\limits_{j=1}^me_j\cS^j\f(\ux)=2\int\limits_\Gamma E_0(\uy-\ux)\underline{n}(\uy)(\mathcal{M}_{\uy}\tilde{f}(\uy))d\uy,
\end{equation}
\begin{equation}\label{eqseg}
\frac{\mu+\lambda}{2}\sum\limits_{j=1}^me_j\cS^j\f(\ux)+\frac{3\mu+\lambda}{2}\sum\limits_{j=1}^m\cS^j\f(\ux)e_j=2\int\limits_\Gamma (\overline{\mathcal{M}}_{\uy}\tilde{f}(\uy))\underline{n}(\uy)E_0(\uy-\ux)d\uy,
\end{equation}
which come from similar analysis to that in the proof of \eqref{relutil} and \eqref{relutil2}. For the sake of completeness, we shall include only the basic idea of the proof of \eqref{eqprim}. 

We have
\begin{align*}
&\frac{\mu+\lambda}{2}\sum\limits_{j=1}^m\cS^j\f(\ux)e_j+\frac{3\mu+\lambda}{2}\sum\limits_{j=1}^me_j\cS^j\f(\ux)\\
&=\frac{\mu+\lambda}{2}f^j(\ux)e_j+\frac{3\mu+\lambda}{2}e_jf^j(\ux)+2\int\limits_\Gamma E_0(\uy-\ux)\underline{n}(\uy)\left(\frac{\mu+\lambda}{2}R_{\ux}(\uy)\puy+\frac{3\mu+\lambda}{2}\puy R_{\ux}(\uy)\right)d\uy.
\end{align*}
We then conclude from
\begin{align*}
&\frac{\mu+\lambda}{2}R_{\ux}(\uy)\puy+\frac{3\mu+\lambda}{2}\puy R_{\ux}(\uy)\\&=\frac{\mu+\lambda}{2}\left(\sum\limits_{j=1}^m f^j(\uy)e_j-\sum\limits_{j=1}^m f^j(\ux)e_j\right)+\frac{3\mu+\lambda}{2}\left(\sum\limits_{j=1}^m e_jf^j(\uy)-\sum\limits_{j=1}^m e_jf^j(\ux)\right)
\end{align*}
and 
$$\int\limits_{\Gamma}E_0(\uy-\ux)\underline{n}(\uy) d\uy=\frac{1}{2}$$
that \eqref{eqprim} is proved.

In the next section we derive an important property of $\cS$. Precisely, it is an involution operator on $\Li(1+\alpha,\Gamma)$. This will allow us to appreciate the relations \eqref{eqprim} and \eqref{eqseg} as
\begin{equation}\label{Sl2}
(\mathcal{M}_{\ux}\tilde{\mathcal{S}}_{\mathcal{L}})(\ux)=2\int\limits_\Gamma E_0(\uy-\ux)\underline{n}(\uy)(\mathcal{M}_{\uy}\tilde{f}(\uy))d\uy,\quad \ux\in\Gamma,
\end{equation}
\begin{equation}\label{Sl22}
(\overline{\mathcal{M}}_{\ux}\tilde{\mathcal{S}}_{\mathcal{L}})(\ux)=2\int\limits_\Gamma (\overline{\mathcal{M}}_{\uy}\tilde{f}(\uy))\underline{n}(\uy)E_0(\uy-\ux)d\uy,\quad \ux\in\Gamma,
\end{equation}
where, in this case, $\tilde{\mathcal{S}}_{\mathcal{L}}$ is the Whitney extension of $\cS\f$. 

\section{Involution property of $\cS$}
An involution operator is a linear operator that is its own inverse. A careful look at specific situation of \eqref{SI} reveals that it may seem  not to be even an injection. Indeed, for instance, the first component $\cS^0$ is defined by $f^{0}$ and a special linear combination of $f^{j}$ that form the $\mathcal{M}$- and $\overline{\mathcal{M}}$-operator of $\tilde{f}$ restricted to $\Gamma$, but not by the whole collection $\f$ itself. 

Addressing this issue a question to ask is whether the values of $f^{0}$ and $\mathcal{M}_{\uy}\tilde{f}(\uy)$ (or $\overline{\mathcal{M}}_{\uy}\tilde{f}(\uy)$) on $\Gamma$ will determine the collection $\f=\{f^{j},\,0\le j\le m\}$ in $\Li(1+\alpha,\Gamma)$. This question is answered positively in the following theorem, which at the same time, is a  rather generalization of \cite[Theorem 3]{ABMM} and plays a crucial role in the proof of our main achievements.   
\begin{theorem}\label{main}
Let $\f\in\Li(1+\alpha,\Gamma)$ be such that $f^{0}=0$ and 
\begin{equation}\label{Df=Dg}
c_1\sum_{j=1}^m f^{j}(\uy)e_j+c_2\sum_{j=1}^m e_jf^{j}(\uy)=0,\,\,\;c_1,c_2\in\R,\; c_1\not=\pm c_2,\;\uy\in\Gamma.
\end{equation}
Then
\[\f=0,\,\mbox{i.e.,}\,\,f^{j}=0,\,\,0\le j\le m.\]
In particular, if $c_1=\frac{\mu+\lambda}{2}$ and $c_2=\frac{3\mu+\lambda}{2}$ then $f^0$ and $\mathcal{M}\tilde{f}$ on $\Gamma$ will determine the function $\f$. 
\end{theorem} 
\pf Before proving our theorem, we have to note that its conclusion does not follow from \eqref{Df=Dg} and the linear independence of the $e_j$'s, unless $\f$  is assumed to be a $\R$-valued function. For general $\R_{0,m}$-valued functions the above argument loses its validity and we need a quite different method.

Let $\ux,\,\uy\in\Gamma$, then there exists a constant $M$ being independent of $\ux,\uy$ such that
\begin{equation}\label{as}
\bigg\|\sum_{j=1}^mf^j(\uy)(x_j-y_j)\bigg\|\le C|\ux-\uy|^{1+\alpha},
\end{equation}
where we have used the assumption $f^0=0$ in $\Gamma$.

Since 
\[
x_j-y_j=-\frac{1}{2}[e_j(\ux-\uy)+(\ux-\uy)e_j],
\]
inequality (\ref{as}) reflects that
\begin{equation}\label{ineq1}
\bigg\|\frac{1}{2}\sum_{j=1}^mf^j(\uy)[e_j(\ux-\uy)+(\ux-\uy)e_j]\bigg\|\le C|\ux-\uy|^{1+\alpha}.
\end{equation}
%by the assumption $\sum_{j=1}^m[g^j(\uy)-f^{j}(\uy)]e_j=0$ in $\Gamma$.

On dividing both sides of \eqref{ineq1} through by $|\ux-\uy|$, we see that

\begin{equation}\label{ineq2}
\bigg\|\frac{1}{2}\sum_{j=1}^m f^j(\uy)\left[e_j\frac{(\ux-\uy)}{|\ux-\uy|}+\frac{(\ux-\uy)}{|\ux-\uy|}e_j\right]\bigg\|\le C|\ux-\uy|^{\alpha}.
\end{equation}
The tangent space of $\Gamma$ at $\uy$, denoted by ${\mbox{Tan}}(\Gamma,\uy)$, is generated by the vectors $\displaystyle\frac{(\ux-\uy)}{|\ux-\uy|}$ and letting $\ux\to\uy$ along any arbitrary-chosen trajectory in $\Gamma$.

Let $\{\underline{v}_1,\underline{v}_2,...,\underline{v}_{m-1}\}$ be an orthonormal basis of ${\mbox{Tan}}(\Gamma,\uy)$ and suppose it admits the following canonical representation in the standard basis of $\R^m$:
\begin{align*}
\underline{v}_1&=\alpha_{11}e_1+\alpha_{12}e_2+...+\alpha_{1m}e_m,\\
\underline{v}_2&=\alpha_{21}e_1+\alpha_{22}e_2+...+\alpha_{2m}e_m,\\
\underline{v}_3&=\alpha_{31}e_1+\alpha_{32}e_2+...+\alpha_{3m}e_m,\\
\vdots &\quad\vdots\\
\underline{v}_{m-1}&=\alpha_{(m-1)1}e_1+\alpha_{(m-1)2}e_2+...+\alpha_{(m-1)m}e_m.
\end{align*} 

\noindent 
From \eqref{ineq2} it follows that
 \begin{equation}\label{=0}
\frac{1}{2}\sum_{j=1}^m f^j(\uy)(e_j\underline{v}_k+\underline{v}_ke_j)=0,
\end{equation}
for $k=1,\dots ,m-1$. On the other hand, we have
\begin{align*}
&\frac{1}{2}\sum_{j=1}^m f^j(\underline{y})(e_j\underline{v}_k+\underline{v}_ke_j)=-\sum_{j=1}^mf^j(\uy)\alpha_{kj}.
\end{align*}
\noindent
In this way we arrive to the system
$$\left\{\begin{array}{cc}
\alpha_{11}f^1(\underline{y})+\cdots +\alpha_{1m}f^m(\underline{y})&=0,\\
\alpha_{21}f^1(\underline{y})+\cdots +\alpha_{2m}f^m(\underline{y})&=0,\\
\vdots&\vdots\\
\alpha_{(m-1)1}f^1(\underline{y})+\cdots +\alpha_{(m-1)m}f^m(\underline{y})&=0,
\end{array}\right.$$
\noindent
with associated matrix
\begin{equation}\label{matr}
\mathcal{A}:=\begin{pmatrix}
\alpha_{11}&\alpha_{12}&...&\alpha_{1m}\\
\alpha_{21}&\alpha_{22}&...&\alpha_{2m}\\
\vdots&\vdots&\vdots&\vdots\\
\alpha_{(m-1)1}&\alpha_{(m-1)2}&...&\alpha_{(m-1)m}
\end{pmatrix}.
\end{equation}
Since the range of $\mathcal{A}$ is $m-1$, specific column of (\ref{matr}), say the $2^{th}$ column, can be expressed as a linear combination of its $m-1$ remaining columns.

So we have
\begin{align*}
\alpha_{12}&=\gamma_1\alpha_{11}+\gamma_2\alpha_{13}+...+\gamma_{m-1}\alpha_{1m},\\
\alpha_{22}&=\gamma_1\alpha_{21}+\gamma_2\alpha_{23}+...+\gamma_{m-1}\alpha_{2m},\\
\alpha_{32}&=\gamma_1\alpha_{31}+\gamma_2\alpha_{33}+...+\gamma_{m-1}\alpha_{3m},\\
\vdots&\quad\vdots\\
\alpha_{(m-1)2}&=\gamma_1\alpha_{(m-1)1}+\gamma_2\alpha_{(m-1)3}+...+\gamma_{m-1}\alpha_{(m-1)m},
\end{align*}
and the above system becomes
\begin{align*}
\alpha_{11}[f^1(\underline{y})+\gamma_1f^2(\underline{y})]+\alpha_{13}[f^3(\underline{y})+\gamma_2f^2(\underline{y})]+\cdots+\alpha_{1m}[f^m(\underline{y})+\gamma_{m-1}f^2(\underline{y})]&=0,\\
\alpha_{21}[f^1(\underline{y})+\gamma_1f^2(\underline{y})]+\alpha_{23}[f^3(\underline{y})+\gamma_2f^2(\underline{y})]+\cdots+\alpha_{2m}[f^m(\underline{y})+\gamma_{m-1}f^2(\underline{y})]&=0,\\
\alpha_{31}[f^1(\underline{y})+\gamma_1f^2(\underline{y})]+\alpha_{33}[f^3(\underline{y})+\gamma_2f^2(\underline{y})]+\cdots+\alpha_{3m}[f^m(\underline{y})+\gamma_{m-1}f^2(\underline{y})]&=0,\\
\vdots&\quad\vdots\\
\alpha_{(m-1)1}[f^1(\underline{y})+\gamma_1f^2(\underline{y})]+\alpha_{(m-1)3}[f^3(\underline{y})+\gamma_2f^2(\underline{y})]+\cdots+\alpha_{(m-1)m}[f^m(\underline{y})+\gamma_{m-1}f^2(\underline{y})]&=0.
\end{align*} 

\noindent
As the range of the associated $(m-1)\times(m-1)$ matrix to the above homogeneous system is $m-1$, only the trivial solution exists. Consequently
\begin{equation}\label{igualdades}
\left\{\begin{array}{ccc}
f^1(\underline{y})&=-\gamma_1f^2(\underline{y}),\\
f^3(\underline{y})&=-\gamma_2f^2(\underline{y}),\\
\vdots&\quad\vdots\\
f^m(\underline{y})&=-\gamma_{m-1}f^2(\underline{y}).
\end{array}\right.
\end{equation}
Then 
\begin{align*}
&c_1f^2(\underline{y})(-\gamma_1e_1+e_2-\gamma_2e_3+...+\gamma_{m-1}e_m)+c_2(-\gamma_1e_1+e_2-\gamma_2e_3+...+\gamma_{m-1}e_m)f^2(\uy)\\&=c_1\sum_{j=1}^mf^j(\underline{y})e_j+c_2\sum_{j=1}^me_jf^j(\uy)=0,
\end{align*}
under assumption. Let us denote vector $-\gamma_1e_1+e_2-\gamma_2e_3+...+\gamma_{m-1}e_m$ by $\underline{u}$. Since $\underline{u}$ is a non-vanishing invertible Clifford vector, whose multiplicative inverse is $-\underline{u}/|\underline{u}|^2$, then we have the following system 
\begin{equation}\label{syste24}
\left\{\begin{array}{rl}
c_1f^2(\uy)-\frac{c_2}{|\underline{u}|} \underline{u}f^2(\uy)\underline{u}=0,\\
c_2f^2(\uy)-\frac{c_1}{|\underline{u}|}\underline{u}f^2(\uy)\underline{u}=0.
\end{array}\right.
\end{equation}
Multiplying the first and second equations of \eqref{syste24} by $c_1$ and $c_2$, respectively, and subtracting, we obtain that $(c_1^2-c_2^2)f^2(\uy)=0$, hence $f^2(\uy)=0$. Finally, the remaining equalities are a direct consequence of \eqref{igualdades}.\qed 
\begin{remark}
Note that if $c_1=\pm c_2$ in \eqref{Df=Dg} then this relation can degenerate and not imply the determination of $\f$ on $\Gamma$. For example, if $\f$ is an $\R$-valued function and $c_1=-c_2$, then relation \eqref{Df=Dg} is satisfied automatically and does not determine the  $f^j$'s in the collection. Whereas if $m$ is even, $c_1=c_2$, and $\f$ is a collection of pseudoscalar fields, then we also arrive at the same conclusion as before.
\end{remark}

\subsection{Plemelj-Privalov Theorem}
The invariance of the H\"older classes under the action of the complex Cauchy singular integral operator is established by the well-known Plemelj-Privalov theorem, named after the Slovenian mathematician Josip Plemelj and the Russian Ivan
Ivanovich Privalov. The result was obtained by Plemelj \cite{plemelj} for the case
of smooth curves and rediscovered by Privalov \cite{privalov} for the circle and, subsequently,
for every piecewise-smooth curve without cusps \cite{privalov2}. This theorem is crucial in the holomorphic function theory, especially for solving boundary value problems. It provides an important tool for analyzing singular operators on the boundary of regions, which has applications in various areas of physics and engineering. Another remarkable work in this area is that of Dyn'kin \cite{dynkin}, which deals with estimates of the smoothness of Cauchy integral in terms of moduli of  smoothness of the function. Due to the variety of uses of this topic, new attempts were made to generalize the Plemelj-Privalov theorem to higher dimensions.

We now derive a Plemelj-Privalov-type theorem for $\cS$, which is of interest by itself. The proof runs as that of \cite[Theorem 2]{DAB} or \cite[Theorem 4]{ABMM} so only we sketch the most part of it.
\begin{theorem}\label{PP}
The operator $\cS$ keeps invariant the Lipschitz class $\Li(1+\alpha,\Gamma)$, i.e.,
\begin{equation}\label{inclusion}
\cS (\Li(1+\alpha,\Gamma))\subset\Li(1+\alpha,\Gamma).
\end{equation}
\end{theorem}
\pf The task is to prove that the collection $\bh:=\{h^j:\,j=0,\dots,m\}$, with
\begin{align}\label{h0}
&h^0(\ux)=2\int\limits_\Gamma E_0(\uy-\ux)\underline{n}(\uy)R_{\ux}(\uy)d\uy-\frac{3\mu+\lambda}{\mu(2\mu+\lambda)}\int\limits_\Gamma E_1(\uy-\ux)\underline{n}(\uy)[\mathcal{M}_{\uy}R_{\ux}(\uy)]d\uy\nonumber\\
&-\frac{\mu+\lambda}{2\mu(2\mu+\lambda)}\left(\int\limits_\Gamma E_0(\uy-\ux)[\overline{\mathcal{M}}_{\uy}R_{\ux}(\uy)]\underline{n}(\uy)(\uy-\ux)d\uy+\sum_{j=1}^me_j\int\limits_\Gamma E_1(\uy-\ux)[\overline{\mathcal{M}}_{\uy}R_{\ux}(\uy)]\underline{n}(\uy)d\uy e_j\right),
\end{align}
and 
\begin{align}\label{hj}
&h^{j}(\ux)=2\int\limits_\Gamma E_0^j(\uy-\ux)\underline{n}(\uy)R_{\ux}(\uy)d\uy-\frac{3\mu+\lambda}{\mu(2\mu+\lambda)}\int\limits_\Gamma E_1^j(\uy-\ux)\underline{n}(\uy)[\mathcal{M}_{\uy}R_{\ux}(\uy)]d\uy\nonumber\\
&-\frac{\mu+\lambda}{2\mu(2\mu+\lambda)}\left(\int\limits_\Gamma E_0^j(\uy-\ux)[\overline{\mathcal{M}}_{\uy}R_{\ux}(\uy)]\underline{n}(\uy)(\uy-\ux)d\uy+\sum_{j=1}^me_j\int\limits_\Gamma E_1^j(\uy-\ux)[\overline{\mathcal{M}}_{\uy}R_{\ux}(\uy)]\underline{n}(\uy)d\uy e_j\right)\nonumber\\
&+\frac{\mu+\lambda}{2\mu(2\mu+\lambda)}\int_\Gamma E_0(\uy-\ux)[\overline{\mathcal{M}}_{\uy}R_{\ux}(\uy)]\underline{n}(\uy)e_jd\uy
\end{align}
for $j=1,\dots,m$, belongs to $\Li(1+\alpha,\Gamma)$.

Let $\ut,\up\in\Gamma$ such that $|\ut-\up|=2r$. To prove \eqref{inclusion} we proceed to show that
\begin{equation}\label{h1}
\|h^0(\ut)-h^0(\up)-\sum\limits_{j=1}^m h^j(\up)(t_j-p_j)\|\le M|\ut-\up|^{1+\alpha}
\end{equation}
and
\begin{equation}\label{h2}
\|h^j(\ut)-h^j(\up)\|\le M|\ut-\up|^\alpha,\,j=1,\dots,m
\end{equation}
for a constant $M$ to be independent of $\ut,\up$.

Let us split $\Gamma$ into three disjoint pieces: $\Gamma_1=\Gamma\cap B_r(\ut)$, $\Gamma_2=\Gamma\cap B_r(\up)$ and $\Gamma_3=\Gamma\setminus\Gamma_1\cup\Gamma_2$ (see Fig. \ref{figure}). To simplify notation we let $h^j_i$ stand for the analogous term  to $h^j$ in \eqref{h0}-\eqref{hj}, where the domain of integration is restricted to $\Gamma_i$. It is sufficient to show that conditions \eqref{h1} and \eqref{h2} hold for each $h^j_i$ with $i=1,2,3$. A brief outline of the proof for \eqref{h1} is given below, \eqref{h2} may be handled in much the same way.        
\begin{figure}[h!]\label{figure}
\centering
\includegraphics[scale=1]{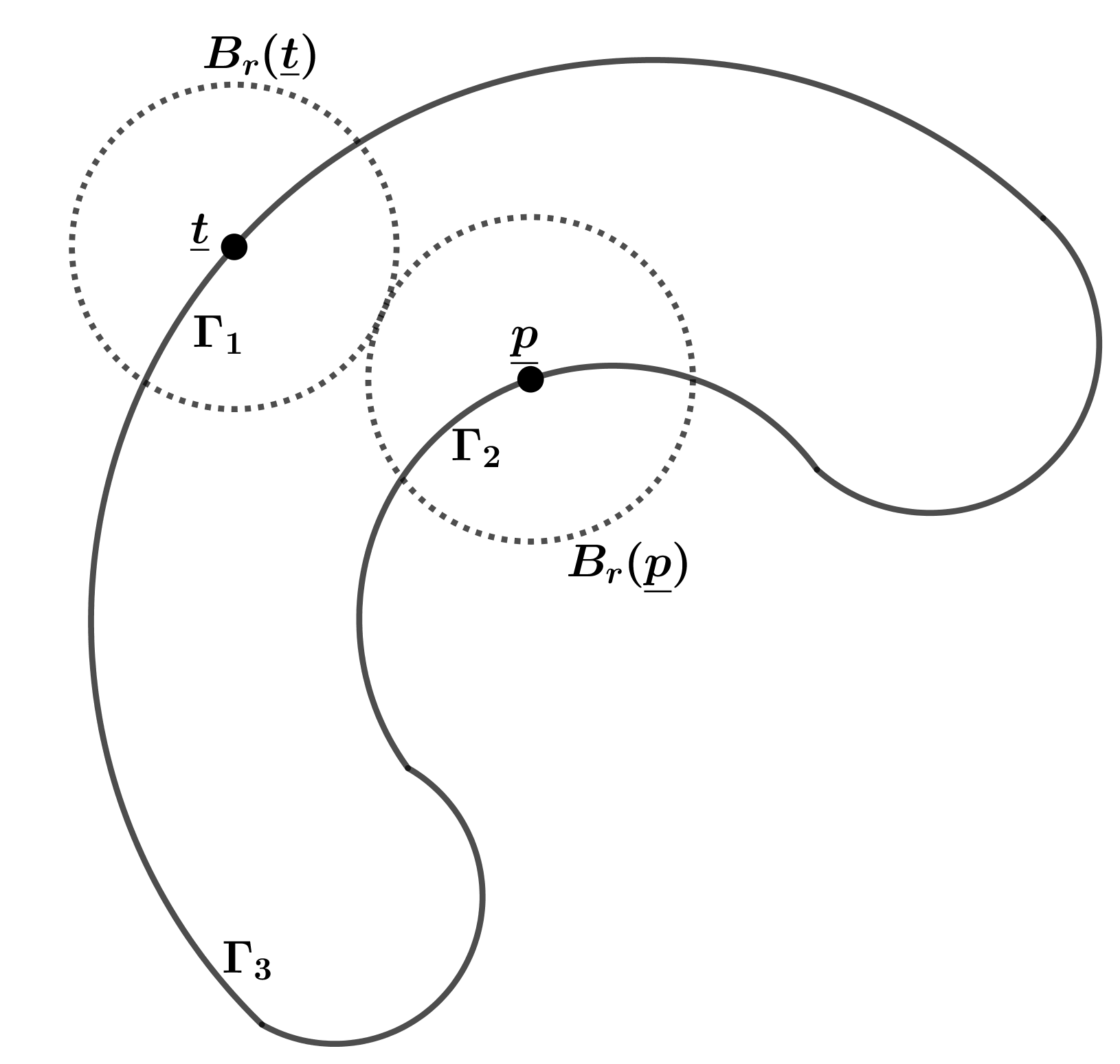}
\caption{Surface decomposition}
\end{figure}

We need only consider the case $i=3$; the other cases $i=1,2$ follow very closely the proof of \cite[Theorem 2]{DAB}. To this end, observe that
\[
R_{\ut}(\uy)=R_{\up}(\uy)+R_{\ut}(\up)+\sum\limits_{j=1}^m[f^j(\up)-f^j(\ut)](y_j-p_j)
\]
and 
\[
R_{\ut}(\uy)\puy=\tilde{f}(\uy)\puy-\tilde{f}(\ut)\piut,\quad \puy R_{\ut}(\uy)=\puy\tilde{f}(\uy)-\piut\tilde{f}(\ut),
\]
and, therefore,
$$\mathcal{M}_{\uy}R_{\ut}(\uy)=\mathcal{M}_{\uy}\tilde{f}(\uy)-\mathcal{M}_{\ut}\tilde{f}(\ut),$$
$$\overline{\mathcal{M}}_{\uy}R_{\ut}(\uy)=\overline{\mathcal{M}}_{\uy}\tilde{f}(\uy)-\overline{\mathcal{M}}_{\ut}\tilde{f}(\ut).$$

A short computation shows that
\[
h^0_3(\ut)-h^0_3(\up)-\sum\limits_{j=1}^m h^j_3(\up)(t_j-p_j)=I_1+I_2,
\]
where
\begin{align*}\label{I1}
&I_1=2\int\limits_{\Gamma_3} \left[E_0(\uy-\ut)-E_0(\uy-\up)-\sum_{j=1}^mE_0^j(\uy-\up)(t_j-p_j)\right]\underline{n}(\uy)R_{\up}(\uy)d\uy\\&-\frac{3\mu+\lambda}{\mu(2\mu+\lambda)}\int\limits_{\Gamma_3} \left[E_1(\uy-\ut)-E_1(\uy-\up)-\sum_{j=1}^mE_1^j(\uy-\up)(t_j-p_j)\right]\underline{n}(\uy)[\mathcal{M}_{\uy}R_{\up}]d\uy\nonumber\\
&-\frac{\mu+\lambda}{2\mu(2\mu+\lambda)}\int\limits_{\Gamma_3}\left[E_0(\uy-\ut)(\overline{\mathcal{M}}_{\uy}R_{\up})n(\uy)(\uy-\ut)-E_0(\uy-\up)(\overline{\mathcal{M}}_{\uy}R_{\up})n(\uy)(\uy-\up)\right]d\uy
\\
&+\frac{\mu+\lambda}{2\mu(2\mu+\lambda)}\left(\int\limits_\Gamma\sum_{j=1}^m\partial_{z_j}[ E_0(\uy-\uz)(\overline{\mathcal{M}}_{\uy}R_{\up})\underline{n}(\uy)(\uy-\uz)]|_{\uz=\up}(t_j-p_j)d\uy\right)\nonumber\\&-\frac{\mu+\lambda}{2\mu(2\mu+\lambda)}\left[\sum_{i=1}^me_i\int\limits_{\Gamma_3}\left(E_1(\uy-\ut)-E_1(\uy-\up)-\sum_{j=1}^mE_1^j(\uy-\up)(t_j-p_j)\right)(\overline{\mathcal{M}}_{\uy}R_{\up})\underline{n}(\uy)d\uy e_i\right]\nonumber
\end{align*}
and 
\begin{eqnarray*}\label{I2}
&I_2=2\int\limits_{\Gamma_3} E_0(\uy-\ut)\underline{n}(\uy)\left[R_{\ut}(\up)+\sum\limits_{j=1}^m(f^j(\up)-f^j(\ut))(y_j-p_j)\right]d\uy\\&-\frac{3\mu+\lambda}{\mu(2\mu+\lambda)}\int\limits_{\Gamma_3} E_1(\uy-\ut)\underline{n}(\uy)[\mathcal{M}_{\up}R_{\ut}]d\uy\nonumber\\
&-\frac{\mu+\lambda}{2\mu(2\mu+\lambda)}\left(\int\limits_{\Gamma_3} E_0(\uy-\ut)[\overline{\mathcal{M}}_{\up}R_{\ut}]\underline{n}(\uy)(\uy-\ut)d\uy+\sum_{j=1}^me_j\int\limits_{\Gamma_3} E_1(\uy-\ut)[\overline{\mathcal{M}}_{\up}R_{\ut}]\underline{n}(\uy)d\uy e_j\right).
\end{eqnarray*}
To estimate $\|I_2\|$ we make use of the simple observation that  
$$g(\uy)=R_{\ut}(\up)+\sum\limits_{j=1}^m[f^j(\up)-f^j(\ut)](y_j-p_j)$$ 
is a Lam\'e-Navier solution such that 
\[\mathcal{M}_{\uy}g(\uy)=\mathcal{M}_{\up}R_{\ut},\quad \overline{\mathcal{M}}_{\uy}g(\uy)=\overline{\mathcal{M}}_{\up}R_{\ut}.\] 
Then by \eqref{CLame} we can use the standard argument of passage from the integration on $\Gamma_3$ to those on pieces of the spheres $\partial B_r(\ut)$, $\partial B_r(\up)$ lying in $\Omega$.

For $\uy\in\partial B_r(\ut)$ the Clifford norm of each integral in \eqref{I2} is dominated by $|\ut-\up|^{1+\alpha}$, which follows immediately after passing the norm under the integral sign. Moreover, for $\uy\in\partial B_r(\up)$ we have
\[
r\le|\uy-\ut|\le |\uy-\up|+|\up-\ut|=3r
\]    
and the same argument can be used, so $\|I_2\|$ is certainly dominated by $|\ut-\up|^{1+\alpha}$.

Similar considerations apply to $\|I_1\|$ like that employed in \cite[Theorem 2]{DAB}. In our context the crucial point is that functions $$E_0(\uy-\uz), E_0(\uy-\uz)(\overline{\mathcal{M}}_{\uy}R_{\up})\underline{n}(\uy)(\uy-\uz),E_1(\uy-\uz)$$ for $\uy\in\Gamma_3$ are (as a function of $\uz$) infinitely differentiable in a suitable neighborhood of $\ut$. This allows us to apply the Taylor formula, yielding an appropriate bound for the remaining terms and finally the desired norm estimate for each integrals in $I_1$. \qed   

\subsection{$\cS$ is a linear involution}
\begin{theorem}
The operator $\cS:\emph{\Li}(1+\alpha,\Gamma)\mapsto \emph{\Li}(1+\alpha,\Gamma)$ is a linear involution. That is,
\begin{equation}\label{invo}
[\cS^2\f]^{j}=\f^{j}\,\,\mbox{for all}\,\, j=0,\dots, m.
\end{equation}
\end{theorem}
\pf In accordance with Theorem \ref{main}, the procedure is to prove that for $\ux\in\Gamma$
\begin{itemize}
\item[(i)]\[[\cS^2\f]^{0}(\ux)=f^{0}(\ux),\]
\item[(ii)] \[\frac{\mu+\lambda}{2}\sum_{j=1}^m [\cS^2\f]^{j}(\ux)e_j+\frac{3\mu+\lambda}{2}\sum_{j=1}^m e_j[\cS^2\f]^{j}(\ux)=\frac{\mu+\lambda}{2}\sum_{j=1}^m f^{j}(\ux)e_j+\frac{3\mu+\lambda}{2}\sum_{j=1}^m e_jf^{j}(\ux).
\]
\end{itemize}
Let us first show the following equality
\begin{align*}
&[\cS^2\f]^{0}(\ux)=2\int\limits_\Gamma E_0(\uy-\ux)\underline{n}(\uy)\cS^0\f(\uy)d\uy-\frac{3\mu+\lambda}{\mu(2\mu+\lambda)}\int\limits_\Gamma E_1(\uy-\ux)\underline{n}(\uy)[\mathcal{M}_{\uy}\tilde{\cS}(\uy)]d\uy\nonumber\\
&-\frac{\mu+\lambda}{2\mu(2\mu+\lambda)}\left(\int\limits_\Gamma E_0(\uy-\ux)[\overline{\mathcal{M}}_{\uy}\tilde{\cS}(\uy)]\underline{n}(\uy)(\uy-\ux)d\uy+\sum_{j=1}^me_j\int\limits_\Gamma E_1(\uy-\ux)[\overline{\mathcal{M}}_{\uy}\tilde{\cS}(\uy)]\underline{n}(\uy)d\uy e_j\right).
\end{align*}
Hence, by \eqref{eqprim} and \eqref{eqseg} we obtain

\begin{align*}
[\cS^2\f]^{0}(\ux)&=2\int\limits_\Gamma 2\int\limits_\Gamma E_0(\uy-\ux)\underline{n}(\uy)E_0(\uz-\uy)\underline{n}(\uz)f^0(\uz)d\uz d\uy\\&\quad-\frac{2(3\mu+\lambda)}{\mu(2\mu+\lambda)}\int\limits_\Gamma\int\limits_\Gamma E_0(\uy-\ux)\underline{n}(\uy)E_1(\uz-\uy)\underline{n}(\uz)[\mathcal{M}_{\uz}\tilde{f}(\uz)]d\uz d\uy\\
&\quad-\frac{\mu+\lambda}{\mu(2\mu+\lambda)}\int\limits_\Gamma\int\limits_\Gamma E_0(\uy-\ux)\underline{n}(\uy)E_0(\uz-\uy)[\overline{\mathcal{M}}_{\uz}\tilde{f}(\uz)]\underline{n}(\uz)(\uz-\uy)d\uz d\uy\\
&\quad-\frac{\mu+\lambda}{\mu(2\mu+\lambda)}\sum_{j=1}^m\int\limits_\Gamma\int\limits_\Gamma E_0(\uy-\ux)\underline{n}(\uy)e_jE_1(\uz-\uy)[\overline{\mathcal{M}}_{\uz}\tilde{f}(\uz)]\underline{n}(\uz)e_jd\uz d\uy\\
&\quad-\frac{2(3\mu+\lambda)}{\mu(2\mu+\lambda)}\int\limits_\Gamma\int\limits_\Gamma E_1(\uy-\ux)\underline{n}(\uy)E_0(\uz-\uy)\underline{n}(\uz)[\mathcal{M}_{\uz}\tilde{f}(\uz)]d\uz d\uy\\
&\quad-\frac{\mu+\lambda}{\mu(2\mu+\lambda)}\int\limits_\Gamma\int\limits_\Gamma E_0(\uy-\ux)[\overline{\mathcal{M}}_{\uz}\tilde{f}(\uz)]\underline{n}(\uz)E_0(\uz-\uy)\underline{n}(\uy)(\uy-\ux)d\uz d\uy\\
&\quad-\frac{\mu+\lambda}{\mu(2\mu+\lambda)}\sum_{j=1}^m\int\limits_\Gamma\int\limits_\Gamma e_jE_1(\uy-\ux)[\overline{\mathcal{M}}_{\uz}\tilde{f}(\uz)]\underline{n}(\uz)E_0(\uz-\uy)\underline{n}(\uy)e_jd\uz d\uy.
\end{align*}
Since 
\[
2\int\limits_\Gamma 2\int\limits_\Gamma E_0(\uy-\ux)\underline{n}(\uy)E_0(\uz-\uy)\underline{n}(\uz)f^0(\uz)d\uz d\uy=f^0(\ux),
\]
it will thus be sufficient to prove that 
\begin{align*}
I&:=-\frac{2(3\mu+\lambda)}{\mu(2\mu+\lambda)}\int\limits_\Gamma\int\limits_\Gamma E_0(\uy-\ux)\underline{n}(\uy)E_1(\uz-\uy)\underline{n}(\uz)[\mathcal{M}_{\uz}\tilde{f}(\uz)]d\uz d\uy\\
&\quad-\frac{\mu+\lambda}{\mu(2\mu+\lambda)}\int\limits_\Gamma\int\limits_\Gamma E_0(\uy-\ux)\underline{n}(\uy)E_0(\uz-\uy)[\overline{\mathcal{M}}_{\uz}\tilde{f}(\uz)]\underline{n}(\uz)(\uz-\uy)d\uz d\uy\\
&\quad-\frac{\mu+\lambda}{\mu(2\mu+\lambda)}\sum_{j=1}^m\int\limits_\Gamma\int\limits_\Gamma E_0(\uy-\ux)\underline{n}(\uy)e_jE_1(\uz-\uy)[\overline{\mathcal{M}}_{\uz}\tilde{f}(\uz)]\underline{n}(\uz)e_jd\uz d\uy\\
&\quad-\frac{2(3\mu+\lambda)}{\mu(2\mu+\lambda)}\int\limits_\Gamma\int\limits_\Gamma E_1(\uy-\ux)\underline{n}(\uy)E_0(\uz-\uy)\underline{n}(\uz)[\mathcal{M}_{\uz}\tilde{f}(\uz)]d\uz d\uy\\
&\quad-\frac{\mu+\lambda}{\mu(2\mu+\lambda)}\int\limits_\Gamma\int\limits_\Gamma E_0(\uy-\ux)[\overline{\mathcal{M}}_{\uz}\tilde{f}(\uz)]\underline{n}(\uz)E_0(\uz-\uy)\underline{n}(\uy)(\uy-\ux)d\uz d\uy\\
&\quad-\frac{\mu+\lambda}{\mu(2\mu+\lambda)}\sum_{j=1}^m\int\limits_\Gamma\int\limits_\Gamma e_jE_1(\uy-\ux)[\overline{\mathcal{M}}_{\uz}\tilde{f}(\uz)]\underline{n}(\uz)E_0(\uz-\uy)\underline{n}(\uy)e_jd\uz d\uy=0.
\end{align*}
A close inspection of the above formula reveals that each integrand there represents the product of a Cauchy kernel with a weakly singular ones. Then we can change of the order of integration, which after a short computation yields
\[
I=\int\limits_{\Gamma} I(\ux,\uz)d\uz,
\]
where
\begin{align*}
I(\ux,\uz)&:=-\frac{2(3\mu+\lambda)}{\mu(2\mu+\lambda)}\int\limits_\Gamma E_0(\uy-\ux)\underline{n}(\uy)E_1(\uz-\uy)\underline{n}(\uz)[\mathcal{M}_{\uz}\tilde{f}(\uz)] d\uy\\
&\quad-\frac{\mu+\lambda}{\mu(2\mu+\lambda)}\int\limits_\Gamma E_0(\uy-\ux)\underline{n}(\uy)E_0(\uz-\uy)[\overline{\mathcal{M}}_{\uz}\tilde{f}(\uz)]\underline{n}(\uz)(\uz-\uy) d\uy\\
&\quad-\frac{\mu+\lambda}{\mu(2\mu+\lambda)}\sum_{j=1}^m\int\limits_\Gamma E_0(\uy-\ux)\underline{n}(\uy)e_jE_1(\uz-\uy)[\overline{\mathcal{M}}_{\uz}\tilde{f}(\uz)]\underline{n}(\uz)e_j d\uy\\
&\quad-\frac{2(3\mu+\lambda)}{\mu(2\mu+\lambda)}\int\limits_\Gamma E_1(\uy-\ux)\underline{n}(\uy)E_0(\uz-\uy)\underline{n}(\uz)[\mathcal{M}_{\uz}\tilde{f}(\uz)] d\uy\\
&\quad-\frac{\mu+\lambda}{\mu(2\mu+\lambda)}\int\limits_\Gamma E_0(\uy-\ux)[\overline{\mathcal{M}}_{\uz}\tilde{f}(\uz)]\underline{n}(\uz)E_0(\uz-\uy)\underline{n}(\uy)(\uy-\ux) d\uy\\
&\quad-\frac{\mu+\lambda}{\mu(2\mu+\lambda)}\sum_{j=1}^m\int\limits_\Gamma e_jE_1(\uy-\ux)[\overline{\mathcal{M}}_{\uz}\tilde{f}(\uz)]\underline{n}(\uz)E_0(\uz-\uy)\underline{n}(\uy)e_j d\uy.
\end{align*}
Let $\epsilon>0$ and consider the balls $B_\epsilon(\ux)$ and $B_\epsilon(\uz)$ of radius $\epsilon$ centered at $\ux$ and $\uz$, respectively. 

Then, with the notation $\Gamma_\epsilon$ for $\Gamma\setminus(B_\epsilon(\ux)\cup B_\epsilon(\uz))$ we have
\begin{align*}
I(\ux,\uz)=\lim_{\epsilon\to 0}I_\epsilon(\ux,\uz):=\lim_{\epsilon\to 0}\bigg[&-\frac{2(3\mu+\lambda)}{\mu(2\mu+\lambda)}\int\limits_{\Gamma_\epsilon} E_0(\uy-\ux)\underline{n}(\uy)E_1(\uz-\uy)\underline{n}(\uz)[\mathcal{M}_{\uz}\tilde{f}(\uz)] d\uy\\
&\quad-\frac{\mu+\lambda}{\mu(2\mu+\lambda)}\int\limits_{\Gamma_\epsilon} E_0(\uy-\ux)\underline{n}(\uy)E_0(\uz-\uy)[\overline{\mathcal{M}}_{\uz}\tilde{f}(\uz)]\underline{n}(\uz)(\uz-\uy) d\uy\\
&\quad-\frac{\mu+\lambda}{\mu(2\mu+\lambda)}\sum_{j=1}^m\int\limits_{\Gamma_\epsilon} E_0(\uy-\ux)\underline{n}(\uy)e_jE_1(\uz-\uy)[\overline{\mathcal{M}}_{\uz}\tilde{f}(\uz)]\underline{n}(\uz)e_j d\uy\\
&\quad-\frac{2(3\mu+\lambda)}{\mu(2\mu+\lambda)}\int\limits_{\Gamma_\epsilon} E_1(\uy-\ux)\underline{n}(\uy)E_0(\uz-\uy)\underline{n}(\uz)[\mathcal{M}_{\uz}\tilde{f}(\uz)] d\uy\\
&\quad-\frac{\mu+\lambda}{\mu(2\mu+\lambda)}\int\limits_{\Gamma_\epsilon} E_0(\uy-\ux)[\overline{\mathcal{M}}_{\uz}\tilde{f}(\uz)]\underline{n}(\uz)E_0(\uz-\uy)\underline{n}(\uy)(\uy-\ux) d\uy\\
&\quad-\frac{\mu+\lambda}{\mu(2\mu+\lambda)}\sum_{j=1}^m\int\limits_{\Gamma_\epsilon} e_jE_1(\uy-\ux)[\overline{\mathcal{M}}_{\uz}\tilde{f}(\uz)]\underline{n}(\uz)E_0(\uz-\uy)\underline{n}(\uy)e_j d\uy\bigg].
\end{align*}
Applying the Stokes formula to $\Omega_\epsilon:=\Omega\setminus(B_\epsilon(\ux)\cup B_\epsilon(\uz))$ and using the equalities
$$\puy[E_0(\uz-\uy)(\overline{\mathcal{M}}_{\uz}\tilde{f}(\uz))\underline{n}(\uz)(\uz-\uy)]=-\sum_{j=1}^me_jE_0(\uz-\uy)(\overline{\mathcal{M}}_{\uz}\tilde{f}(\uz))\underline{n}(\uz)e_j,$$
$$[(\uy-\ux)(\overline{\mathcal{M}}_{\uz}\tilde{f}(\uz))\underline{n}(\uz)E_0(\uz-\uy)]\puy=\sum_{j=1}^me_j (\overline{\mathcal{M}}_{\uz}\tilde{f}(\uz))\underline{n}(\uz)E_0(\uz-\uy)e_j,$$
and 
$$\puy\left[\sum_{j=1}^me_jE_1(\uz-\uy)(\overline{\mathcal{M}}_{\uz}\tilde{f}(\uz))\underline{n}(\uz)e_j\right]=2(\overline{\mathcal{M}}_{\uz}\tilde{f}(\uz))\underline{n}(\uz)E_0(\uz-\uy)+\sum_{j=1}^me_jE_0(\uz-\uy)(\overline{\mathcal{M}}_{\uz}\tilde{f}(\uz))\underline{n}(\uz)e_j,$$
we get
\begin{align*}
&I_\epsilon(\ux,\uz)-\bigg\{-\frac{2(3\mu+\lambda)}{\mu(2\mu+\lambda)}\int\limits_{C_\epsilon(\ux)\cup C_\epsilon(\uz)} E_0(\uy-\ux)\underline{n}(\uy)E_1(\uz-\uy)\underline{n}(\uz)[\mathcal{M}_{\uz}\tilde{f}(\uz)] d\uy\\
&\quad-\frac{\mu+\lambda}{\mu(2\mu+\lambda)}\int\limits_{C_\epsilon(\ux)\cup C_\epsilon(\uz)} E_0(\uy-\ux)\underline{n}(\uy)E_0(\uz-\uy)[\overline{\mathcal{M}}_{\uz}\tilde{f}(\uz)]\underline{n}(\uz)(\uz-\uy) d\uy\\
&\quad-\frac{\mu+\lambda}{\mu(2\mu+\lambda)}\sum_{j=1}^m\int\limits_{C_\epsilon(\ux)\cup C_\epsilon(\uz)} E_0(\uy-\ux)\underline{n}(\uy)e_jE_1(\uz-\uy)[\overline{\mathcal{M}}_{\uz}\tilde{f}(\uz)]\underline{n}(\uz)e_j d\uy\\
&\quad-\frac{2(3\mu+\lambda)}{\mu(2\mu+\lambda)}\int\limits_{C_\epsilon(\ux)\cup C_\epsilon(\uz)} E_1(\uy-\ux)\underline{n}(\uy)E_0(\uz-\uy)\underline{n}(\uz)[\mathcal{M}_{\uz}\tilde{f}(\uz)] d\uy\\
&\quad-\frac{\mu+\lambda}{\mu(2\mu+\lambda)}\int\limits_{C_\epsilon(\ux)\cup C_\epsilon(\uz)} E_0(\uy-\ux)[\overline{\mathcal{M}}_{\uz}\tilde{f}(\uz)]\underline{n}(\uz)E_0(\uz-\uy)\underline{n}(\uy)(\uy-\ux) d\uy\\
&\quad-\frac{\mu+\lambda}{\mu(2\mu+\lambda)}\sum_{j=1}^m\int\limits_{C_\epsilon(\ux)\cup C_\epsilon(\uz)} e_jE_1(\uy-\ux)[\overline{\mathcal{M}}_{\uz}\tilde{f}(\uz)]\underline{n}(\uz)E_0(\uz-\uy)\underline{n}(\uy)e_j d\uy\bigg\}\\
&=\frac{2(3\mu+\lambda)}{\mu(2\mu+\lambda)}\int\limits_{\Omega_\epsilon}E_0(\uy-\ux)E_0(\uz-\uy)\underline{n}(\uz)[\mathcal{M}_{\uz}\tilde{f}(\uz)]d\uy\\&\quad+\frac{\mu+\lambda}{\mu(2\mu+\lambda)}\sum_{j=1}^m\int\limits_{\Omega_\epsilon}E_0(\uy-\ux)e_jE_0(\uz-\uy)[\overline{\mathcal{M}}_{\uz}\tilde{f}(\uz)]\underline{n}(\uz)e_jd\uy\\&\quad-\frac{2(\mu+\lambda)}{\mu(2\mu+\lambda)}\int\limits_{\Omega_\epsilon}E_0(\uy-\ux)[\overline{\mathcal{M}}_{\uz}\tilde{f}(\uz)]\underline{n}(\uz)E_0(\uz-\uy)d\uy\\
&\quad-\frac{\mu+\lambda}{\mu(2\mu+\lambda)}\sum_{j=1}^m\int\limits_{\Omega_\epsilon}E_0(\uy-\ux)e_jE_0(\uz-\uy)[\overline{\mathcal{M}}_{\uz}\tilde{f}(\uz)]\underline{n}(\uz)e_jd\uy\\
&\quad-\frac{2(3\mu+\lambda)}{\mu(2\mu+\lambda)}\int\limits_{\Omega_\epsilon}E_0(\uy-\ux)E_0(\uz-\uy)\underline{n}(\uz)[\mathcal{M}_{\uz}\tilde{f}(\uz)]d\uy\\
&\quad-\frac{\mu+\lambda}{\mu(2\mu+\lambda)}\sum_{j=1}^m\int\limits_{\Omega_\epsilon}e_j[\overline{\mathcal{M}}_{\uz}\tilde{f}(\uz)]\underline{n}(\uz)E_0(\uz-\uy)e_jE_0(\uy-\ux)d\uy\\
&\quad -\frac{\mu+\lambda}{\mu(2\mu+\lambda)}\sum_{j=1}^m\int\limits_{\Omega_\epsilon}e_j[\overline{\mathcal{M}}_{\uz}\tilde{f}(\uz)]\underline{n}(\uz)E_0(\uz-\uy)E_0(\uy-\ux)e_jd\uy,
\end{align*}
where $C_\epsilon(\ux):=\partial{B_\epsilon(\ux)}\cap\Omega$, $C_\epsilon(\uz):=\partial{B_\epsilon(\uz)}\cap\Omega$.

Being aware that $E_0(\uy-\ux)e_i+e_iE_0(\uy-\ux)$ is $\R$-valued, it follows that
\begin{align*}
&-\frac{2(\mu+\lambda)}{\mu(2\mu+\lambda)}\int\limits_{\Omega_\epsilon}E_0(\uy-\ux)[\overline{\mathcal{M}}_{\uz}\tilde{f}(\uz)]\underline{n}(\uz)E_0(\uz-\uy)d\uy\\
&\quad-\frac{\mu+\lambda}{\mu(2\mu+\lambda)}\sum_{j=1}^m\int\limits_{\Omega_\epsilon}e_j[\overline{\mathcal{M}}_{\uz}\tilde{f}(\uz)]\underline{n}(\uz)E_0(\uz-\uy)e_jE_0(\uy-\ux)d\uy\\
&\quad -\frac{\mu+\lambda}{\mu(2\mu+\lambda)}\sum_{j=1}^m\int\limits_{\Omega_\epsilon}e_j[\overline{\mathcal{M}}_{\uz}\tilde{f}(\uz)]\underline{n}(\uz)E_0(\uz-\uy)E_0(\uy-\ux)e_jd\uy\\
&=-\frac{2(\mu+\lambda)}{\mu(2\mu+\lambda)}\int\limits_{\Omega_\epsilon}E_0(\uy-\ux)[\overline{\mathcal{M}}_{\uz}\tilde{f}(\uz)]\underline{n}(\uz)E_0(\uz-\uy)d\uy\\
&\quad-\frac{\mu+\lambda}{\mu(2\mu+\lambda)}\sum_{j=1}^m\int\limits_{\Omega_\epsilon}e_j(e_jE_0(\uy-\ux)+E_0(\uy-\ux)e_j)[\overline{\mathcal{M}}_{\uz}\tilde{f}(\uz)]\underline{n}(\uz)E_0(\uz-\uy)d\uy.
\end{align*}
Now we make use of the identity 
\[
\sum_{i=1}^m e_iE_0(\uy-\ux)e_i=(m-2)E_0(\uy-\ux),
\]
which implies that the last integrals vanish.

Consequently, we have
\begin{align*}
I_\epsilon(\ux,\uz)&=-\frac{2(3\mu+\lambda)}{\mu(2\mu+\lambda)}\int\limits_{C_\epsilon(\ux)\cup C_\epsilon(\uz)} E_0(\uy-\ux)\underline{n}(\uy)E_1(\uz-\uy)\underline{n}(\uz)[\mathcal{M}_{\uz}\tilde{f}(\uz)] d\uy\\
&\quad-\frac{\mu+\lambda}{\mu(2\mu+\lambda)}\int\limits_{C_\epsilon(\ux)\cup C_\epsilon(\uz)} E_0(\uy-\ux)\underline{n}(\uy)E_0(\uz-\uy)[\overline{\mathcal{M}}_{\uz}\tilde{f}(\uz)]\underline{n}(\uz)(\uz-\uy) d\uy\\
&\quad-\frac{\mu+\lambda}{\mu(2\mu+\lambda)}\sum_{j=1}^m\int\limits_{C_\epsilon(\ux)\cup C_\epsilon(\uz)} E_0(\uy-\ux)\underline{n}(\uy)e_jE_1(\uz-\uy)[\overline{\mathcal{M}}_{\uz}\tilde{f}(\uz)]\underline{n}(\uz)e_j d\uy\\
&\quad-\frac{2(3\mu+\lambda)}{\mu(2\mu+\lambda)}\int\limits_{C_\epsilon(\ux)\cup C_\epsilon(\uz)} E_1(\uy-\ux)\underline{n}(\uy)E_0(\uz-\uy)\underline{n}(\uz)[\mathcal{M}_{\uz}\tilde{f}(\uz)] d\uy\\
&\quad-\frac{\mu+\lambda}{\mu(2\mu+\lambda)}\int\limits_{C_\epsilon(\ux)\cup C_\epsilon(\uz)} E_0(\uy-\ux)[\overline{\mathcal{M}}_{\uz}\tilde{f}(\uz)]\underline{n}(\uz)E_0(\uz-\uy)\underline{n}(\uy)(\uy-\ux) d\uy\\
&\quad-\frac{\mu+\lambda}{\mu(2\mu+\lambda)}\sum_{j=1}^m\int\limits_{C_\epsilon(\ux)\cup C_\epsilon(\uz)} e_jE_1(\uy-\ux)[\overline{\mathcal{M}}_{\uz}\tilde{f}(\uz)]\underline{n}(\uz)E_0(\uz-\uy)\underline{n}(\uy)e_j d\uy
\end{align*}
Taking advantage of the fact that for $\uy\in C_\epsilon(\ux)$ ($\uy\in C_\epsilon(\uz)$) the normal vector $\underline{n}(\uy)$ has the explicit form $\displaystyle\frac{\uy-\ux}{|\uy-\ux|}$ ($\displaystyle\frac{\uy-\uz}{|\uy-\uz|}$), we obtain
\begin{align*}
&I(\ux,\uz)=\lim_{\epsilon\to 0}I_\epsilon(\ux,\uz)=-\frac{3\mu+\lambda}{\mu(2\mu+\lambda)}E_1(\uz-\ux)n(\uz)[\mathcal{M}_{\uz}\tilde{f}(\uz)]-\frac{\mu+\lambda}{2\mu(2\mu+\lambda)}E_0(\uz-\ux)[\overline{\mathcal{M}}_{\uz}\tilde{f}(\uz)]n(\uz)(\uz-\ux)\\
&\quad-\frac{\mu+\lambda}{2\mu(2\mu+\lambda)}\sum_{j=1}^me_jE_1(\uz-\ux)[\overline{\mathcal{M}}_{\uz}\tilde{f}(\uz)]n(\uz)e_j+\frac{3\mu+\lambda}{\mu(2\mu+\lambda)}E_1(\uz-\ux)n(\uz)[\mathcal{M}_{\uz}\tilde{f}(\uz)]d\uy\\
&\quad+\frac{\mu+\lambda}{2\mu(2\mu+\lambda)}E_0(\uz-\ux)[\overline{\mathcal{M}}_{\uz}\tilde{f}(\uz)]n(\uz)(\uz-\ux)+\frac{\mu+\lambda}{2\mu(2\mu+\lambda)}\sum_{j=1}^me_jE_1(\uz-\ux)[\overline{\mathcal{M}}_{\uz}\tilde{f}(\uz)]n(\uz)e_j\\
&=0,
\end{align*}
which proves (i). 

Let us now examine (ii), nearly identical to obvious from \eqref{eqprim}. Indeed, we have
\begin{align*}
&\frac{\mu+\lambda}{2}\sum\limits_{j=1}^m[\cS^2\f]^{j}(\ux)e_j+\frac{3\mu+\lambda}{2}\sum\limits_{j=1}^me_j[\cS^2\f]^{j}(\ux)\\&=\frac{\mu+\lambda}{2}\sum\limits_{j=1}^m[\cS^j(\cS\f)](\ux)e_j+\frac{3\mu+\lambda}{2}\sum\limits_{j=1}^me_j[\cS^j(\cS\f)](\ux)\\&=2\int\limits_\Gamma E_0(\uy-\ux)\underline{n}(\uy)[\mathcal{M}_{\uy}\tilde{\cS}(\uy)]d\uy\nonumber\\&=2\int\limits_\Gamma 2\int\limits_\Gamma E_0(\uy-\ux)\underline{n}(\uy)E_0(\uz-\uy)n(\uz)[\mathcal{M}_{\uz}\tilde{f}(\uz)]d\uz d\uy\nonumber\\&=\frac{\mu+\lambda}{2}\sum_{j=1}^m f^{j}(\ux)e_j+\frac{3\mu+\lambda}{2}\sum_{j=1}^m e_jf^{j}(\ux),
\end{align*}
as claimed.\qed
\section{Hardy Decomposition}
On account of the above section, $\cP^+:=\frac{1}{2}(\I+\cS)$ and $\cP^-:=\frac{1}{2}(\I-\cS)$ are projection operators on $\Li(1+\alpha,\Gamma)$. Namely,
\[
\cP^+\cP^+=\cP^+,\,\,\cP^-\cP^-=\cP^-,\,\,\cP^+\cP^-=\cP^-\cP^+=0.
\]
Therefore, the Lipschitz class $\Li(1+\alpha,\Gamma)$ admits the orthogonal decomposition
\[\Li(1+\alpha,\Gamma)=\Li^+(1+\alpha,\Gamma)\oplus\Li^-(1+\alpha,\Gamma),\]
where
\[
\Li^+(1+\alpha,\Gamma):=\mbox{im}\cP^+,\,\Li^-(1+\alpha,\Gamma):=\mbox{im}\cP^-.
\]
The rest of the paper is devoted to characterize the structure of $\Li^\pm(1+\alpha,\Gamma)$.
\begin{theorem}
A collection $\f\in\Li(1+\alpha,\Gamma)$ belongs to $\Li^+(1+\alpha,\Gamma)$ if and only if there exists a Lam\'e-Navier solution $F$ in $\Omega_+$ which together with $\mathcal{M}_{\ux}F$, continuously extend to $\Gamma$ and satisfy
\begin{equation}\label{traceconditions}
F|_\Gamma=f^0,\,\,\mathcal{M}_{\ux}F|_\Gamma=\frac{\mu+\lambda}{2}\sum_{j=1}^mf^j e_j+\frac{3\mu+\lambda}{2}\sum_{j=1}^me_jf^j.
\end{equation}
\end{theorem}
\pf Assume first that $\f\in\Li^+(1+\alpha,\Gamma)$, then there exists $\bg\in\Li(1+\alpha,\Gamma)$ such that 
\begin{equation}\label{+}
\f=\frac{1}{2}(\cS\bg+\bg). 
\end{equation}
Let us introduce the function 
\[
F(\ux)=[\C^0\bg](\ux),\,\,\ux\in\Omega_+.
\]
Obviously, $F$ is a Lam\'e-Navier solution in $\Omega_+$ and for $\uz\in\Gamma$ we have
\[
F(\uz)=[\C^0\bg]^+(\uz)=\frac{1}{2}(\cS^0\bg+g^0)=f^0.
\]
Moreover, by \eqref{i-l} we have:
\[
\mathcal{M}_{\ux}F(\ux)=\mathcal{M}_{\ux}[\C^0\bg](\ux)={\cC}^l\left[\frac{\mu+\lambda}{2}\sum_{j=1}^mg^j e_j+\frac{3\mu+\lambda}{2}\sum_{j=1}^me_jg^j\right](\ux),\,\,\ux\in\Omega_+.
\]
Consequently, for $\uz\in\Gamma$
\begin{align*}
(\mathcal{M}_{\ux}F)(\uz)&=\left[{\cC}^l\left[\frac{\mu+\lambda}{2}\sum_{j=1}^mg^j e_j+\frac{3\mu+\lambda}{2}\sum_{j=1}^me_jg^j \right]\right]^+(\uz)\\&=\frac{1}{2} \left[2\int\limits_{\Gamma}E_0(\uy-\uz)\underline{n}(\uy)\left(\frac{\mu+\lambda}{2}\sum_{j=1}^m g^j(\uy)e_j+\frac{3\mu+\lambda}{2}\sum_{j=1}^m e_jg^j(\uy)\right) d\uy\right.\\&\quad\left.+\frac{\mu+\lambda}{2}\sum_{j=1}^m g^j(\uz)e_j+\frac{3\mu+\lambda}{2}\sum_{j=1}^m e_jg^j(\uz)\right]\\&=\frac{1}{2}\sum\limits_{j=1}^m\left[\frac{\mu+\lambda}{2}\cS^j\bg(\uz)e_j+\frac{3\mu+\lambda}{2}e_j\cS^j\bg(\uz)+\frac{\mu+\lambda}{2}g^j(\uz)e_j+\frac{3\mu+\lambda}{2}e_jg^j(\uz)\right].
\end{align*}
The last equality is established by \eqref{eqprim}.

Finally, we conclude from \eqref{+} that
\[
(\mathcal{M}_{\ux}F)(\uz)=\frac{\mu+\lambda}{2}\sum_{j=1}^mf^j(\uz) e_j+\frac{3\mu+\lambda}{2}\sum_{j=1}^me_jf^j(\uz),\,\uz\in\Gamma.
\]
Now suppose that such Lam\'e-Navier solution $F$ exists and satisfies \eqref{traceconditions}. Then, formula \eqref{CLame} implies $F(\ux)=\C^0\f(\ux)$, $\ux\in\Omega_+$. We claim that $\cP^+\f=\f$ and hence $\f\in\Li^+(1+\alpha,\Gamma)$. Indeed, we have
\[
[\cP^+\f]^0(\uz)=\frac{1}{2}(\cS^0\f(\uz)+f^0(\uz))=\lim_{\Omega_+\ni\ux\to\uz}\C^0\f(\ux)=F(\uz)=f^0(\uz).
\]
On the other hand
\begin{align*}
&\frac{\mu+\lambda}{2}\sum_{j=1}^m[\cP^+\f]^j(\uz)e_j+\frac{3\mu+\lambda}{2}\sum_{j=1}^me_j[\cP^+\f]^j(\uz)\\
&=\frac{1}{2}\left[\frac{\mu+\lambda}{2}\sum_{j=1}^m\cS^j\f(\uz)e_j+\frac{3\mu+\lambda}{2}\sum_{j=1}^me_j\cS^j\f(\uz)+\frac{\mu+\lambda}{2}\sum_{j=1}^m f^j(\uz)e_j+\frac{3\mu+\lambda}{2}\sum_{j=1}^m e_jf^j(\uz)\right]\\&=\frac{1}{2} \left[2\int\limits_{\Gamma}E_0(\uy-\uz)\underline{n}(\uy)\left(\frac{\mu+\lambda}{2}\sum_{j=1}^m f^j(\uy)e_j+\frac{3\mu+\lambda}{2}\sum_{j=1}^m e_jf^j(\uy)\right) d\uy\right.\\&\quad\left.+\frac{\mu+\lambda}{2}\sum_{j=1}^m f^j(\uz)e_j+\frac{3\mu+\lambda}{2}\sum_{j=1}^m e_jf^j(\uz)\right].
\end{align*}
The last equality is a consequence of \eqref{eqprim}.

Under the assumption that $\frac{\mu+\lambda}{2}\sum_{j=1}^m f^j(\uz)e_j+\frac{3\mu+\lambda}{2}\sum_{j=1}^m e_jf^j(\uz)$ is a H\"older continuous function in $\Gamma$, which represents the interior limit value of the monogenic function $\mathcal{M}_{\ux}F$, we have 
\begin{align*}
&\frac{1}{2} \left[2\int\limits_{\Gamma}E_0(\uy-\uz)\underline{n}(\uy)\left(\frac{\mu+\lambda}{2}\sum_{j=1}^m f^j(\uy)e_j+\frac{3\mu+\lambda}{2}\sum_{j=1}^m e_jf^j(\uy)\right) d\uy\right.\\&\quad\left.+\frac{\mu+\lambda}{2}\sum_{j=1}^m f^j(\uz)e_j+\frac{3\mu+\lambda}{2}\sum_{j=1}^m e_jf^j(\uz)\right]\\
&=\frac{\mu+\lambda}{2}\sum_{j=1}^m f^j(\uz)e_j+\frac{3\mu+\lambda}{2}\sum_{j=1}^m e_jf^j(\uz),
\end{align*}
which implies
\[
\frac{\mu+\lambda}{2}\sum_{j=1}^m[\cP^+\f]^j(\uz)e_j+\frac{3\mu+\lambda}{2}\sum_{j=1}^me_j[\cP^+\f]^j(\uz)=\frac{\mu+\lambda}{2}\sum_{j=1}^m f^j(\uz)e_j+\frac{3\mu+\lambda}{2}\sum_{j=1}^m e_jf^j(\uz).
\]
Therefore $[\cP^+\f]^j(\uz)=f^j(\uz)$ by Theorem \ref{main}.\qed

The proof of the following result is analogous and hence it is omitted.
\begin{theorem}
A collection $\f\in\Li(1+\alpha,\Gamma)$ belongs to $\Li^-(1+\alpha,\Gamma)$ if and only if there exists a Lam\'e-Navier  solution $F$ in $\Omega_-$ which together with $\mathcal{M}_{\ux}F$, continuously extend to $\Gamma$, with $F(\infty)=0$, $\|\pux F(\ux)\|=o(\frac{1}{|\ux|})$ as $\ux\to\infty$,  and such that
\[
F|_\Gamma=f^0,\,\,\mathcal{M}_{\ux}F|_\Gamma=\frac{\mu+\lambda}{2}\sum_{j=1}^mf^j e_j+\frac{3\mu+\lambda}{2}\sum_{j=1}^me_jf^j.
\]
\end{theorem}

\section*{Conflict of interest} The authors declare that they have no conflict of interest regarding the publication of this paper.

\section*{Acknowledgements}
Daniel Alfonso Santiesteban gratefully acknowledges the financial support of the Secretar\'ia de Ciencia, Humanidades, Tecnolog\'ia e Innovaci\'on (SECIHTI)  (Grant Number 1043969).

\end{document}